\documentclass{siamltex}
\usepackage{amsmath}
\usepackage{amsfonts}
\usepackage{amssymb}
\usepackage{latexsym}
\usepackage{upgreek}
\usepackage{cleveref}
\usepackage{stmaryrd} 
\usepackage{subfigure}
\usepackage{comment} 
\usepackage{bm}
\usepackage[usenames]{color}
\usepackage[usenames,dvipsnames,svgnames]{xcolor}
\usepackage{graphicx}
\usepackage{amscd}
\usepackage{acronym,arydshln}

\def\u{{\bm u}}

\def\el{\nonumber}
\def\cl{\nonumber \\}

\def\Ijnme{International Journal for Numerical Methods in Engineering}

\def\Siamna{SIAM Journal on Numerical Analysis}

\def\Nm{Numerische Mathematik}

\def\u{\bm u}
\def\vt{\bm v}

\def\u{{\bm u}}

\def\vt{{\bm v}}
\def\n{{\bm n}}

\def\vtres#1{{ \vert\!\vert\!\vert {#1} \vert\!\vert\!\vert }}

\def\eu {{\bm e_{\u} }} 
\def\ceu {{\nabla \times \eu }} 
\def\u{{\bm u}}
\def\vt{{\bm v}}

\def\Ijnme{International Journal for Numerical Methods in Engineering}

\def\Siamna{SIAM Journal on Numerical Analysis}
\def\Nm{Numerische Mathematik}

\def\cl{\nonumber\\}
\def\el{\nonumber}

\def\vtres#1{{\Vert {#1} \Vert}}

\newcommand{\Eq}[1]{(\ref{eq:#1})}

\newcommand\os{\begin{color}{black}}
\newcommand\rc{\begin{color}{black}}
\newcommand\db{\begin{color}{black}}
\newcommand\rcr{\begin{color}{black}}
\newcommand\ot{\end{color}}

\allowdisplaybreaks

\title{Nitsche's prescription of Dirichlet conditions in the finite element approximation of Maxwell's problem}
\author{Daniele Boffi\footnotemark[1]
      \and
       {Ramon Codina}\footnotemark[2]
      \and
        \"Onder T\"urk\footnotemark[3]
}

\begin{document}

\maketitle

\renewcommand{\thefootnote}{\fnsymbol{footnote}}
\footnotetext[1]{King Abdullah University of Science and Technology, Thuwal, Saudi Arabia, and University of Pavia, Pavia, Italy}
\footnotetext[2]{Universitat Polit\`{e}cnica de Catalunya and Centre Internacional de M\`{e}todes Num\`{e}rics en Enginyeria, Barcelona, Spain}
\footnotetext[3]{Middle East Technical University, Ankara, Turkey}
\renewcommand{\thefootnote}{\arabic{footnote}}

\begin{abstract}
In this paper we consider the finite element approximation of Maxwell's problem and analyse the prescription of essential boundary conditions in a weak sense using Nitsche's method. To avoid indefiniteness of the problem, the original equations are augmented with the gradient of a scalar field that allows one to impose the zero divergence of the magnetic induction, even if the exact solution for this scalar field is zero. Two finite element approximations are considered, namely, one in which the approximation spaces are assumed to satisfy the appropriate inf-sup condition that render the standard Galerkin method stable, and another augmented and stabilised one that permits the use of finite element interpolations of arbitrary order. Stability and convergence results are provided for the two finite element formulations considered.
\end{abstract}



\pagestyle{myheadings}
\thispagestyle{plain}
\markboth{D. Boffi, R. Codina and \"O. T\"urk}{Nitsche's prescription of Dirichlet conditions in the FE approximation of Maxwell's problem}

\section{Introduction}

The Maxwell equations govern the electromagnetic wave propagation, and hence are involved in many scientific and industrial fields. A continuous research has been and is being conducted to correctly approximate the solution to the Maxwell problem by means of a number of numerical techniques, among  which the finite element (FE) method is the most widely used. Some studies dealing with FE methods for approximating the solutions to time-harmonic Maxwell problems are \cite{aylwin2023sinum,badia-codina-2009-2,chaumont-frelet2023sinum,du2023,monk-zhang-2020}, and the references therein. 

This paper deals with the FE approximation of the following problem: find a magnetic induction field  $\u:\Omega\longrightarrow\mathbb{R}^d$ and a scalar field $p:\Omega\longrightarrow\mathbb{R}$ solution of the boundary value problem
\begin{alignat}{3}
\nu\nabla\times\nabla\times \u + \nabla p& = {\bm f} \quad && \hbox{in}~\Omega,\label{eq:m1}\\
- \nabla\cdot\u  & = 0\quad && \hbox{in}~\Omega,\label{eq:m2}\\
\n\times\u &= \n\times\bar{\u}\quad && \hbox{on}~\Gamma,\label{eq:m3}\\
p & = {\bar p} := 0 \quad && \hbox{on}~\Gamma,\label{eq:m4}
\end{alignat}
where $\Omega$ is a bounded polyhedral domain of $\mathbb{R}^d$ ($d=2,3$), $\Gamma = \partial \Omega$, $\nu>0$ is a physical parameter (the inverse of the magnetic permeability times the electric conductivity), $\bar{\u}$ is given and $\bm f$ is assumed to be solenoidal.  As usual, the scalar field $p$ is introduced to impose that the FE  approximation to $\u$ be solenoidal, since at the continuous level the solution is $p=0$. We call $p$ the magnetic pseudo-pressure.

\db
To the best of our knowledge, this approach was introduced
in~\cite[Formulation F4]{Kikuchi}, for the study of the Maxwell eigenvalue
problem.
\ot

A possible strategy for  weakly imposing Dirichlet boundary conditions is Nitsche's method (see for example  \cite{juntunen-stenberg-09,stenberg-1995}) which consists in penalising the difference between the unknown and its prescribed value on the boundary, with a proper scaling. The weak form is obtained without assuming that the test functions vanish on $\Gamma$; this leads to a boundary term whose symmetric counterpart is also introduced to preserve the symmetry of the variational formulation of the problem. In the context of discontinuous Galerkin methods, this leads to the well known interior penalty method, usually attributed to the works \cite{arnold-1982-1,O-douglas-dupont,wheeler-1978}. 

\os Even though Nitsche's method is widely applied for interface problems  in electromagnetics (see, e.g., \cite{roppert2020,wang2022} and references)therein , the number of studies applying Nitsche's method to prescribe the Dirichlet boundary conditions in FE approaches for the Maxwell problem is very few in the literature. A form of Nitsche's method is used in \cite{winges2023} where a scattered field formulation of the FE method for Maxwell's equations is considered on a so called Huygens' surface which encloses the scatterer being located in the free space embedding this scatterer. In this reference, the computational domain is partitioned into total-  and scattered-field regions, and  the equivalent electric and magnetic surface currents are incorporated in the weak form by means of Nitsche's method. As for a Nitsche type formulation to  directly handle the boundary conditions in the Maxwell problem, a discrete formulation is proposed in \cite{assous2011}, where the ellipticity 
 of the associated bilinear form in curl-div form (with the inclusion of Nitsche's terms) is shown.  

In this paper, our main interest is to analyse Nitsche's method to prescribe the Dirichlet boundary conditions \Eq{m3} and \Eq{m4} in the FE context, and this is why we have not considered Neumann-type boundary conditions. This prescription of the Dirichlet boundary conditions is done for two FE formulations. In the first one, the Galerkin method is employed and the interpolating spaces for $\u$ and $p$ are assumed to satisfy adequate inf-sup conditions that render the discrete problem stable. These interpolations can be for example N\'ed\'elec's elements for $\u$ and standard nodal continuous interpolations for $p$. For the second formulation, we consider that continuous nodal interpolations are used for both $\u$ and $p$. The Galerkin method in this case is unstable and one has to switch to stabilised FE formulations. The one we consider here was introduced and analysed in~\cite{badia-codina-2009-2}.
\ot 

The second approach, i.e., the use of continuous interpolations for the magnetic induction, is of particular interest. 
\os
It has been well understood since the work of \cite{costabel-dauge-2002-1}
that on domains with re-entrant corners, the standard Galerkin method that is
applied with a subspace of continuous piecewise FE spaces is prone to
producing non-physical solutions \ot
\db if the solenoidal condition is imposed by a penalty method. \ot
\os In other words, there are solutions of the Maxwell equations which cannot be properly approximated with the use of standard conforming FEs. A number of alternatives with different nature has already been introduced to potentially restore the use of standard continuous elements \cite{assous1998,badia-codina-2009-2,costabel-dauge-2002-1}. 

There is also a {\em modelling} difficulty in terms of boundary conditions when approximating Maxwell's problem using continuous nodal elements, and it is related to the fact of preserving conformity. If a node belongs to two edges (respectively faces in 3D) that are not co-aligned (respectively co-planar in 3D), the only way to guarantee that the component of the magnetic field tangent to $\Gamma$ is zero is to prescribe {\em all} the components to zero. On a curved boundary, that would imply to prescribe all the field components to zero at all nodes, unless a ${\mathcal C}^1$ description of the boundary is used. If one defines a `numerical' tangent (typically from a numerical normal) and prescribes the resulting tangent components, conformity will not hold exactly. This sort of {\em variational crime} does not appear using a weak prescription of boundary conditions using the type of techniques presented in the following sections. 
\ot     

The paper is organised as follows. In Section~\ref{sec:twofe} the FE approximation to the problem with exact imposition of boundary conditions is presented. Nothing is new there, one of the methods is the standard Galerkin method and the other one is the formulation proposed in~\cite{badia-codina-2009-2}. Section~\ref{sec:nitsche} presents the application of Nitsche's method in combination with the formulations introduced in Section~\ref{sec:twofe}. \os The analysis here is achieved by a novel strategy of approximate imposition of boundary conditions which consists of splitting the continuous functions into one function that vanishes on the boundary and its complement, an idea originally introduced in \cite{codina-baiges-2008-1}.  \ot  Apart from the presented analysis, the novelty is the combination of Nitsche's method with the stabilised formulation introduced in Section~\ref{sec:twofe}, which regarding the stabilisation mechanism is similar but not identical to the discontinuous G
 alerkin formulation presented in~\cite{perugias1}. Numerical results are presented in Section~\ref{sec:numex}, and finally conclusions are drawn in Section~\ref{sec:conclu}.

\section{Two finite element approximations for Maxwell's problem}
\label{sec:twofe}

\subsection{Continuous problem}

Let us introduce some notation. If $X$ is a Hilbert space of functions defined on $\Omega$ where the unknown is sought, its norm is denoted as $\Vert\cdot\Vert_X$, its dual as $X'$, and the duality by $\langle\cdot,\cdot \rangle_\Omega$. If $\Lambda$ is the space of traces on $\Gamma$ of functions in $X$ and $\Lambda'$ is its dual, the duality in this case is written as $\langle\cdot,\cdot \rangle_\Gamma$. The $L^2$-inner product in a domain $\omega$ is denoted by $(\cdot,\cdot )_\omega$. The $L^2(\Omega)$-projection onto a space $X$ is written as $P_X$. Moreover, inequalities up to {\em dimensionless} constants, independent also of the discretisation, are written as $\lesssim$ and $\gtrsim$ for $\leq$ and $\geq$, respectively.

The differential operator of Maxwell's differential equations \Eq{m1}-\Eq{m2} can be written as ${\mathcal L}([\u,p]) = [\nu\nabla\times\nabla\times \u +  \nabla p , - \nabla\cdot\u]$, and then those equations are ${\mathcal L}([\u,p]) = [{\bm f}, 0]$. Let $\vt$ and $q$ be arbitrary functions with the same regularity as $\u$ and $p$, respectively. For future use, after appropriate integration by parts and assuming enough regularity of the functions involved, we get the identity
\begin{align}
\langle {\mathcal L}([\u,p]) , [\vt , q]\rangle_\Omega
=  B([\u,p] , [\vt , q]) - \langle {\mathcal F}_n ( [\u , p]) , {\mathcal D} ([\vt , q])\rangle_\Gamma,\label{eq:byparts}
\end{align}
where
\begin{align}
& B([\u,p],[\vt,q]) = \nu (\nabla\times \u, \nabla\times\vt)_{\Omega} + (\nabla p,\vt)_{\Omega} 
+ (\nabla q,\u)_{\Omega}\label{eq:defB}\\
& {\mathcal F}_n  ([\u , p]) =  [\nu P_t(\nabla\times\u), \n\cdot\u],
\quad {\mathcal F} ([\u,p]) = [\nu \nabla\times\u , \u],\label{eq:defF} \\
& {\mathcal D} ([\vt , q]) = [\n\times \vt, q],\label{eq:defD}
\end{align}
and we have introduced the tangent projection $P_t$ on the boundary $\Gamma$, defined for any vector field $\bm a$ as $P_t({\bm a}) = {\bm a} -({\bm a}\cdot\n) \n$. We could also have defined  ${\mathcal F}_n  ([\u , p] )=  [\nu\, \n\times\nabla\times\u, \n\cdot\u]$, ${\mathcal D} ([\vt , q]) = [P_t(\vt), q]$; the expression chosen is due to the boundary conditions \Eq{m3}-\Eq{m4} that we wish to impose. The need for introducing $P_t$ is merely technical. In the functional spaces where the problem is well posed (see below), only the tangent component of $\nabla\times \u$ is well defined, in the sense that it belongs to the dual space of the trace of ${\bm n}\times \vt$ on $\Gamma$. However, since $P_t(\nabla\times\u) \cdot (\n\times \vt) = (\nabla\times\u) \cdot ( \n\times \vt)$, we omit the projection $P_t$ in what follows.

The variational form of problem \Eq{m1}-\Eq{m4} is well posed in the space $X = V\times Q := H({\rm curl}; \Omega)\times H^1(\Omega)$, where $H({\rm curl}; \Omega)$ is the space of vector fields in $L^2(\Omega)^d$ with curl in $L^2(\Omega)^d$. The subspace made of vectors $\vt\in H({\rm curl}; \Omega)$ such that $\n\times\vt = {\bf 0}$ on $\Gamma$ is denoted by $V_0 = H_0({\rm curl}; \Omega)$, and the subspace of scalar functions in $H^1(\Omega)$ vanishing on $\Gamma$ as $Q_0 = H_0^1(\Omega)$. The space of traces is $\Lambda = H^{-1/2}({\rm div}_\Gamma;\Gamma)\times H^{1/2}(\Gamma)$, the trace operator being $[\u,p]\mapsto [\n\times\u,p]$; for a characterisation of $H^{-1/2}({\rm div}_\Gamma;\Gamma)$ for polyhedral domains and different results about traces of $V$, see e.g. \cite{buffa-et-al-2002} and references therein. 

To ease the writing of the continuous problem, let us consider for the moment $\bar{\u} = {\bf 0}$. The general case could be treated using the standard lifting of $\bar{\u}$ to a function defined on the whole $\Omega$. The weak form of problem \Eq{m1}-\Eq{m4}, imposing the boundary conditions in an essential manner, reads as follows: find $[\u,p]\in X_0 := V_0\times Q_0$ such that 
\begin{align}
& B([\u,p],[\vt,q]) = \langle {\bm f} , \vt\rangle_\Omega,\label{eq:probcont}
\end{align}
for all $[\vt,q]\in X_0$, i.e., $[\vt,q]\in X$ and ${\mathcal D} ([\vt , q]) =
[{\bf 0}, 0]$. This problem is known to be well posed, \db in particular, \ot it holds:
\begin{align}
\inf_{[\u,p]\in X_0}\sup_{[\vt,q]\in X_0} \frac{B([\u,p],[\vt,q])} {\Vert [\u,p] \Vert_X \Vert [\vt,q] \Vert_X} \geq K_B > 0,\label{eq:infsupB}
\end{align}
where $K_B$ is a positive constant and
\begin{align}
& \Vert [\u,p] \Vert^2_X := \Vert \u \Vert^2_V + \Vert p \Vert^2_Q,\cl
& \Vert \u \Vert^2_V  := \nu \Vert\nabla\times \u \Vert^2_{L^2(\Omega)} + \frac{\nu}{L_0^2} \Vert \u \Vert^2_{L^2(\Omega)}, \cl
& \Vert p \Vert^2_Q := \frac{L_0^2}{\nu} \Vert \nabla p \Vert^2_{L^2(\Omega)}, \el
\end{align}
where $L_0$ is a characteristic length of $\Omega$. Note that $\Vert \cdot\Vert_V$ is the norm in $H({\rm curl}; \Omega)$ with adequate scaling coefficients and $\Vert \cdot\Vert_Q$ is a scaled norm in $H_0^1(\Omega)$ because of the Poincar\'e-Friedrichs inequality. In all what follows, scaling coefficients are introduced to make all terms dimensionally consistent.

The continuous problem \Eq{probcont} is equivalent to the two variational equations:
\begin{alignat}{3}
a(\u,\vt) + b(p,\vt) & = \langle {\bm f} , \vt\rangle_\Omega \quad &&\forall \vt\in V_0,\label{eq:varc1}\\
b(q,\u) & = 0\quad  && \forall q\in Q_0,\label{eq:varc2}
\end{alignat}
with 
\begin{align*}
a(\u,\vt) := \nu (\nabla\times \u, \nabla\times\vt)_{\Omega}, \quad  b(p,\vt)  := (\nabla p,\vt)_{\Omega}.
\end{align*}
The inf-sup condition \Eq{infsupB} is then a consequence of the `little' inf-sup condition
\begin{align}
\inf_{p\in Q_0}\sup_{\vt\in V_0} \frac{b(p,\vt)} {\Vert p \Vert_Q \Vert \vt \Vert_V} \geq K_b > 0,\label{eq:littleinfsup}
\end{align}
and the coercivity of $a(\u,\vt)$ in $K_V = \{\vt\in V_0 ~\vert~    b(q,\vt) = 0 ~ \forall q\in Q_0\}$.

\subsection{Galerkin finite element approximation}

Let us consider now the Ga\-ler\-kin FE approximation of problem \Eq{probcont}.
For that, let us construct a FE  partition of $\Omega$, ${\mathcal T}_h = \{
K\}$, with $h = \max_K \{ h_K = {\rm diam} (K), K \in {\mathcal T}_h\}$, which
we consider \rcr shape regular. \ot We assume that the domain $\Omega$
is polyhedral, and that $\Omega = {\rm int}(\bigcup_{K\in {\mathcal T}_h} K)$
for all $h>0$, \os where ${\rm int}(\omega)$ is interior of $\omega$. \ot From ${\mathcal T}_h$ we may construct now FE spaces $V_{h,0} \subset V_0$ and $Q_{h,0} \subset Q_0$, i.e., we consider conforming FE approximations. For any FE function $v_h$, piecewise polynomial of degree $k$, the following inverse and trace inequalities hold:
\begin{align}
& \Vert \nabla v_h \Vert_{L^2(K)} \leq C_{\rm inv} \frac{k^2}{h_K} \Vert v_h\Vert_{L^2(K)}, \label{eq:invineq}\\
& \Vert v_h \Vert_{L^2(\partial K)} \leq C_{\rm trace} \frac{k}{h^{1/2}_K} \Vert v_h\Vert_{L^2(K)}, \label{eq:traceineq}
\end{align}
\rcr
with the corresponding obvious global counterparts:
\begin{align*}
\Vert \nabla v_h \Vert_{L^2(\Omega)} \leq C_{\rm inv} \sum_K \frac{k^2}{h_K} \Vert v_h\Vert_{L^2(K)}, \qquad
\Vert v_h \Vert_{L^2(\Gamma)} \leq C_{\rm trace} \sum_K \frac{k}{h_K^{1/2}} \Vert v_h\Vert_{L^2(K)}.
\end{align*}
To avoid overloading the notation, we shall use the shortcut
\begin{align*}
\sum_K f({h_K}) \Vert v_h\Vert_{L^2(K)}\equiv f({h}) \Vert v_h\Vert_{L^2(\Omega)},
\end{align*}
for any function $f$, understanding that when $f(h)$ is multiplying a global norm, it should be replaced by $f(h_K)$ multiplying the same norm on each $K$ and summing for $K \in {\mathcal T}_h$. The same comment holds when the factor of $h$ multiplies an inner product or a norm of functions defined on $\Gamma$. Note, however, that our results do not assume that the FE partition is quasi-uniform.
\ot

Since we take $k$ fixed in the following analysis, we may consider it absorbed in the constants $C_{\rm inv}$ and $C_{\rm trace}$. We will also make use of the inverse inequality:
\begin{align}
& \Vert v_h \Vert_{L^\infty(\partial K)} \leq  {C_{\rm inv}}{h^{-(d-1)/2}_K} \Vert v_h\Vert_{L^2(\partial K)}. \label{eq:invineq-b}
\end{align}

We assume in this subsection that spaces $V_{h,0}$ and $Q_{h,0}$ satisfy the discrete version of condition \Eq{littleinfsup}, i.e.,
\begin{align}
\inf_{p_h\in Q_{h,0}}\sup_{\vt_h\in V_{h,0}} \frac{b(p_h,\vt_h)} {\Vert p_h \Vert_Q \Vert \vt_h \Vert_V} \geq K_b > 0,\el
\end{align}
or, equivalently, 
\begin{align}
\forall p_h \in Q_{h,0} ~\exists\, \vt_h\in V_{h,0} ~\hbox{such that}~ b(p_h,\vt_h) \geq K_b \Vert p_h \Vert_Q \Vert \vt_h \Vert_V.\label{eq:littleinfsup-h}
\end{align}

As already mentioned, examples of pairs of spaces satisfying this condition
are those based on N\'ed\'elec's elements to construct $V_{h,0}$ and nodal
Lagrangian continuous elements to construct $Q_{h,0}$. In general, condition
\Eq{littleinfsup-h} \db is guaranteed if the diagram\ot
\vskip0.2cm
\begin{center}
$\begin{CD}
 H^1_0(\Omega) @> \nabla >> H_0({\rm curl};\Omega) \\
@VV P_{Q_h} V @VV P_{V_h} V  \\
Q_{h,0} @> \nabla >> V_{h,0}  \\
\end{CD}$
\end{center}
\vskip0.2cm
\noindent
is commutative (see~\cite{douglas-roberts,bbf,hiptmair-2002}).

As for the continuous problem, condition \Eq{littleinfsup-h} and the coercivity of $a(\u,\vt)$ in the discrete version of the kernel $K_V$ also imply the discrete counterpart of \Eq{infsupB}, which we may write as
\begin{align}
& \forall [\u_h,p_h] \in V_{h,0}\times Q_{h,0} ~\exists\, [\vt_{h,0} , q_{h,0}] \in V_{h,0} \times Q_{h,0}  ~\hbox{such that}~ \cl
& \qquad B([\u_h,p_h],[\vt_{h,0} , q_{h,0}]) \geq K_B \Vert [\u_h,p_h] \Vert_{V\times Q}  \Vert [ \vt_{h,0} , q_{h,0} ] \Vert_{V\times Q}.\label{eq:biginfsup-h}
\end{align}

If this inf-sup condition holds, the following problem is well posed: find $[\u_h,p_h]\in V_{h,0} \times Q_{h,0}$ such that
\begin{alignat}{3}
a(\u_h,\vt_h) + b(p_h,\vt_h) & = \langle {\bm f} , \vt_h\rangle_\Omega \quad &&\forall \vt_h\in V_{h,0},\label{eq:gal1}\\
b(q_h,\u_h) & = 0\quad  && \forall q_h\in Q_{h,0}.\label{eq:gal2}
\end{alignat}
This problem admits a unique solution that depends continuously on the data ${\bm f}$. The exact solution for $p_h$ is $p_h = 0$, but introducing it allows one to eliminate the indefiniteness associated to the curl-curl operator. 

We have the following stability and convergence result~\cite{boffi2003,monk-zhang-2020}:

\begin{theorem}\label{th:th1} Suppose that both $V_{h,0}$ and $Q_{h,0}$ are an
inf-sup stable pair satisfying condition \Eq{biginfsup-h}. Then, problem \Eq{gal1}-\Eq{gal2} is well posed, in the sense that it admits a unique solution $[\u_h,p_h]\in V_{h,0}\times Q_{h,0}$ that satisfies
\begin{align*}
\Vert [\u_h,p_h] \Vert_{V\times Q} \lesssim \Vert {\bm f}\Vert_{V'}.
\end{align*}
Furthermore, $[\u_h,p_h]$ converges optimally as $h \to 0$ to the solution $[\u,p]\in V_0 \times Q_0$ of the continuous problem \Eq{varc1}-\Eq{varc2}, in the following sense:
\begin{align}
& \Vert [\u - \u_h, p - p_h] \Vert_{V\times Q}
\lesssim  \inf_{[ \vt_h, q_h] \in V_{h,0}\times Q_{h,0}} \Vert [\u - \vt_h, p - q_h] \Vert_{V\times Q}.\label{eq:errest1}
\end{align}
\end{theorem}

\subsection{Stabilised FE approximation}\label{sec:stab-exactBC}

An alternative to using inf-sup stable spaces $V_{h,0}$-$Q_{h,0}$ is to use a stabilised FE formulation. In this case the situation is somehow particular, as since the solution for the magnetic pseudo-pressure is $p = 0$, there is a `stabilisation' term that can in fact be introduced at the continuous level. Indeed, the solution to problem \Eq{m1}-\Eq{m4} is the same as the solution to 
\begin{alignat}{3}
\nu\nabla\times\nabla\times \u + \nabla p& = {\bm f} \quad && \hbox{in}~\Omega,\label{eq:ma1}\\
\frac{L_0^2}{\nu}\Delta p  - \nabla\cdot\u  & = 0\quad && \hbox{in}~\Omega,\label{eq:ma2}\\
\n\times\u &= \n\times\bar{\u}\quad && \hbox{on}~\Gamma,\label{eq:ma3}\\
p & = {\bar p} := 0 \quad && \hbox{on}~\Gamma,\label{eq:ma4}
\end{alignat}
where the term $-\frac{L_0^2}{\nu}\Delta p$ helps to stabilise the pressure gradient; this problem can be considered an {\em augmented} version of \Eq{m1}-\Eq{m4}. However, when the FE approximation is considered it is also necessary to stabilise the divergence of the discrete velocity. The final stabilised FE formulation we consider in this paper was introduced and analysed in \cite{badia-codina-2009-2} and it consists of finding $[\u_h,p_h]\in V_{h,0} \times Q_{h,0}$ such that
\begin{alignat}{3}
a(\u_h,\vt_h) + b(p_h,\vt_h) + s_u(\u_h,\vt_h) & = \langle {\bm f} , \vt_h\rangle_\Omega \quad &&\forall \vt_h\in V_{h,0},\label{eq:aug1}\\
b(q_h,\u_h) + s_p (p_h,q_h)& = 0\quad  && \forall q_h\in Q_{h,0},\label{eq:aug2}
\end{alignat}
where 
\begin{align}
s_u(\u_h,\vt_h) := c_u  \frac{\nu h^2}{L_0^2}(\nabla\cdot\u_h  , \nabla\cdot \vt_h)_\Omega,\quad
s_p(p_h,q_h) := -\frac{L_0^2}{\nu} (\nabla p_h , \nabla q_h)_\Omega,\label{eq:stabterms}
\end{align}
$c_u$ being an algorithmic constant. \rcr Recall that $h^2 (\nabla\cdot\u_h  , \nabla\cdot \vt_h)_\Omega$ has to be understood as $ \sum_K h_K^2(\nabla\cdot\u_h  , \nabla\cdot \vt_h)_K$.\ot 

Let us write the bilinear form that defines the problem as
\begin{align}
B_{\rm S} ([\u_h,p_h], [\vt_h,q_h])
:= B([\u_h,p_h], [\vt_h,q_h]) + s_u(\u_h,\vt_h)  + s_p (p_h,q_h).\label{eq:defBS}
\end{align}
In \cite{badia-codina-2009-2} it was proved directly that problem \Eq{aug1}-\Eq{aug2} is stable, without the need of proving an inf-sup condition for $B_{\rm S}$. However, this inf-sup condition will be very convenient in the forthcoming analysis, and therefore we prove it here:

\begin{theorem}\label{th:inf-sup-BS}
The bilinear form $B_{\rm S}$ is inf-sup stable in $V_{h,0} \times Q_{h,0}$ in the norm $\Vert \cdot \Vert_{V\times Q,{\rm S}}$ defined as
\begin{align}
\Vert  [\vt_h,q_h] \Vert_{V\times Q,{\rm S}}^2  = \Vert  [\vt_h,q_h] \Vert_{V\times Q}^2  + \nu \frac{h^2}{L_0^2} \Vert \nabla\cdot \u_h \Vert_{L^2(\Omega)}^2,\el
\end{align}
that is, for each $[\u_h,p_h]\in V_{h,0} \times Q_{h,0}$ there exists $[\vt_h,q_h]\in V_{h,0} \times Q_{h,0}$ such that
\begin{align}
B_{\rm S} ([\u_h,p_h], [\vt_h,q_h]) \gtrsim \Vert  [\u_h,p_h] \Vert_{V\times Q,{\rm S}}\Vert  [\vt_h,q_h] \Vert_{V\times Q,{\rm S}}.\el
\end{align}
\end{theorem}

\begin{proof}
Let us start noting that 
\begin{align}
B_{\rm S} ([\u_h,p_h], [\u_h,-p_h]) = \nu \Vert \nabla\times \u_h\Vert_{L^2(\Omega)}^2 +  c_u \nu \frac{h^2}{L_0^2} \Vert \nabla\cdot \u_h \Vert_{L^2(\Omega)}^2
+ \frac{L_0^2}{\nu} \Vert \nabla p_h\Vert_{L^2(\Omega)}^2.\label{eq:bs-t1}
\end{align}
It only remains to obtain control on the $L^2(\Omega)$-norm of $\u_h$. For
that, let us consider its Helmholtz decomposition at continuous level
\begin{align}
& \u_h = \u_d + \frac{L^2_0}{\nu}\nabla r, \quad \hbox{with}\quad \nabla\cdot\u_d = 0,\cl
& \n \times \u_d = \n\times \u_h = {\bf 0}, \quad r = 0\quad \hbox{on}~\Gamma.\el
\end{align}
Note that, in general, $\u_d\not\in V_{h,0}$ and $r\not\in Q_{h,0}$.

Since $\n \times \u_d = {\bf 0}$ on $\Gamma$ and $\nabla\cdot\u_d = 0$, from the Poincar\'e-Friedrichs-type inequality $\Vert \nabla\times \u_d \Vert_{L^2(\Omega)} \gtrsim L_0^{-1} \Vert  \u_d \Vert_{L^2(\Omega)}$ and the fact that $\nabla \times \u_d = \nabla \times \u_h$, \Eq{bs-t1} in fact implies that
\begin{align}
B_{\rm S} ([\u_h,p_h], [\u_h,-p_h]) & \gtrsim \nu \Vert \nabla\times \u_h\Vert_{L^2(\Omega)}^2 
+ \frac{\nu}{L_0^2} \Vert \u_d\Vert_{L^2(\Omega)}^2 \cl
& + \nu \frac{h^2}{L_0^2} \Vert \nabla\cdot \u_h \Vert_{L^2(\Omega)}^2
+ \frac{L_0^2}{\nu} \Vert \nabla p_h\Vert_{L^2(\Omega)}^2.\label{eq:bs-t2}
\end{align}
Thus, only the ${L^2(\Omega)}$ control on $\nabla r$ is needed.  

Let $\pi_h(r)\in Q_{h,0}$ be an interpolant of order at least one of $r\in Q_0$. Since we require $\pi_h(r) = 0$ on $\Gamma$, the Scott-Zhang interpolant can be used.  We now have that:
\begin{align}
& B_{\rm S} ([\u_h,p_h], [{\bf 0} ,  \pi_h(r) ]) = - \frac{L^2_0}{\nu}  ( \nabla p_h , \nabla \pi_h(r))_\Omega + (\u_h ,  \nabla \pi_h(r))_\Omega.\label{eq:bs-t3}
\end{align}
Using the $H^1(\Omega)$-stability of the interpolant and Young's inequality:
\begin{align}
 ( \nabla p_h , \nabla \pi_h(r))_\Omega \lesssim \Vert \nabla p_h \Vert_{L^2(\Omega)}  \Vert \nabla r\Vert_{L^2(\Omega)}
 \leq \frac{1}{2\alpha_1} \Vert \nabla p_h \Vert_{L^2(\Omega)}^2 + \frac{\alpha_1}{2} \Vert \nabla r\Vert_{L^2(\Omega)}^2.\el
\end{align}
The second term in \Eq{bs-t3} can be treated as follows:
\begin{align}
(\u_h ,  \nabla \pi_h(r))_\Omega 
& = \Bigl(\u_d + \frac{L^2_0}{\nu}\nabla r , \nabla r \Bigr)_\Omega +  (\u_h,  \nabla \pi_h(r) - \nabla r)_\Omega\cl
& = \frac{L^2_0}{\nu} \Vert \nabla r \Vert_{L^2(\Omega)}^2 - (\nabla\cdot \u_h, \pi_h(r) -  r )_\Omega\cl
& \gtrsim \frac{L^2_0}{\nu} \Vert \nabla r \Vert_{L^2(\Omega)}^2 - \Vert \nabla\cdot \u_h \Vert_{L^2(\Omega)} h \Vert \nabla r\Vert_{L^2(\Omega)}\cl
& \geq  \frac{L^2_0}{\nu} \Vert \nabla r \Vert_{L^2(\Omega)}^2 - \frac{1}{2\alpha_2} \frac{\nu h^2}{L_0^2} \Vert \nabla\cdot \u_h \Vert_{L^2(\Omega)}^2 - \frac{\alpha_2}{2}  \frac{L^2_0}{\nu} \Vert \nabla r \Vert_{L^2(\Omega)}^2.\el
\end{align}
In the first step we have used the Helmholtz decomposition of $\u_h$, in the second step that $\u_d$ is divergence free and we have integrated by parts the second term, in the third step the approximation property assumed for the interpolant and in the last step Young's inequality.

Using the last bounds in \Eq{bs-t3} and taking $\alpha_1$ and $\alpha_2$ sufficiently small, it turns out that there exists a constant $\gamma \geq 0 $ such that
\begin{align}
& B_{\rm S} ([\u_h,p_h], [{\bf 0} ,  \pi_h(r) ]) \gtrsim \frac{L^2_0}{\nu} \Vert \nabla r \Vert_{L^2(\Omega)}^2 
- \gamma \Bigl(  \frac{L^2_0}{\nu} \Vert \nabla p_h \Vert_{L^2(\Omega)}^2 + \frac{\nu h^2}{L_0^2} \Vert \nabla\cdot \u_h \Vert_{L^2(\Omega)}^2 \Bigr).\el
\end{align}
If we now take $[\vt_h , q_h] = [\u_h , -p_h + \delta  \pi_h(r)]$, with $\delta > 0$ sufficiently small, it follows from this last inequality and from \Eq{bs-t2} that:
\begin{align}
B_{\rm S} ([\u_h,p_h], [\vt_h, q_h]) & \gtrsim \nu \Vert \nabla\times \u_h\Vert_{L^2(\Omega)}^2 
+ \frac{\nu}{L_0^2} \Vert \u_d\Vert_{L^2(\Omega)}^2 
+  \frac{L^2_0}{\nu} \Vert \nabla r \Vert_{L^2(\Omega)}^2 \cl
& + \nu \frac{h^2}{L_0^2} \Vert \nabla\cdot \u_h \Vert_{L^2(\Omega)}^2
+ \frac{L_0^2}{\nu} \Vert \nabla p_h\Vert_{L^2(\Omega)}^2.\el
\end{align}
The $L^2(\Omega)$-orthogonality of the Helmholtz decomposition yields
\begin{align}
 \frac{\nu}{L_0^2} \Vert \u_d\Vert_{L^2(\Omega)}^2 +  \frac{L^2_0}{\nu} \Vert \nabla r \Vert_{L^2(\Omega)}^2 =  \frac{\nu}{L_0^2} \Vert \u_h\Vert_{L^2(\Omega)}^2,\el
\end{align}
and therefore $B_{\rm S} ([\u_h,p_h], [\vt_h, q_h]) \gtrsim \Vert [\u_h , p_h]\Vert_{V\times Q,{\rm S}}^2$. The proof concludes checking that $\Vert [\u_h , p_h]\Vert_{V\times Q,{\rm S}} \gtrsim \Vert [\vt_h , q_h]\Vert_{V\times Q,{\rm S}}$, which again is a consequence of the $H^1(\Omega)$-stability of the interpolant.
\end{proof}

The following results are directly proved in \cite{badia-codina-2009-2}, without using Theorem~\ref{th:inf-sup-BS}:

\begin{theorem}\label{th:th2} Suppose that both $V_{h,0}$ and $Q_{h,0}$ are constructed using continuous nodal based interpolations of arbitrary degree each. Then, problem \Eq{aug1}-\Eq{aug2} is well posed, in the sense that it admits a unique solution $[\u_h,p_h]\in V_{h,0}\times Q_{h,0}$ that satisfies
\begin{align*}
\Vert [\u_h,p_h] \Vert_{V\times Q} \lesssim \Vert {\bm f}\Vert_{V'}.
\end{align*}
Furthermore, $[\u_h,p_h]$ converges optimally as $h \to 0$ to the solution $[\u,p]\in V_0 \times Q_0$ of the continuous problem \Eq{varc1}-\Eq{varc2}, in the following sense:
\begin{align}
& \Vert [\u - \u_h, p - p_h] \Vert_{V\times Q, {\rm S}} \lesssim  \inf_{[ \vt_h, q_h] \in V_{h,0}\times Q_{h,0}} \Vert [\u - \vt_h, p - q_h] \Vert_{V\times Q, {\rm S}}.\label{eq:errest2}
\end{align}
\end{theorem}

The error estimates \Eq{errest1} and \Eq{errest2} are clearly optimal for smooth solutions. In the case of solutions with Sobolev regularity $0 < r < 1$, they are also optimal if the FE meshes are able to interpolate optimally scalar functions of Sobolev regularity $r+1$, whose gradients are components of $\u$. This happens for example if the FE meshes are of Powell-Sabin type (see \cite{badia-codina-2009-2} and references therein for further discussion).

\section{Nitsche's method for Maxwell's problem}
\label{sec:nitsche}

\os In this section we consider \ot that both boundary conditions \Eq{m3} and \Eq{m4} are prescribed weakly, without incorporating them in the FE spaces. Obviously, we may take $\bar{\u} \not = {\bf 0}$, in general, since now assuming homogeneous boundary conditions does not introduce any simplification in the notation.

Here we concentrate on Nitsche's method, which is well understood. It is applied first with the inf-sup stable formulation and later with the stabilised one described earlier. We view Nitsche's method as a way to prescribe boundary conditions, but it can also be understood as a way to impose continuity along interior interfaces, as in the discontinuous Galerkin (dG) method. This is done in particular in \cite{perugias1}, where a dG method is introduced and analysed for the time harmonic Maxwell problem. In fact, in this reference a similar div-div stabilising term as the one in \Eq{stabterms} is employed, the pressure stabilisation being different \os (see also \cite{houston2005})\ot. We also introduce a symmetric form of Nitsche's method, which leads to the so called symmetric interior penalty in the context of dG methods. Finally, since in our case the solution we wish to find is $p = 0$, perhaps it is easier in all cases to prescribe the magnetic pseudo-pressure strongly, but we e
 mploy also Nitsche's strategy in this case, to unify its treatment with that of the magnetic induction $\u$. 

From now on, subscript zero is used for spaces that incorporate homogeneous boundary conditions, whereas this subscript is dropped if no boundary values are prescribed.

\subsection{Nitsche's method using the Galerkin FE approximation} \label{sec:gal-nit}

Let $V_h\subset V$ and $Q_h \subset Q$ be conforming FE spaces, such that the subspaces $V_{h,0}\subset V_h$ and $Q_{h,0}\subset Q_h$ satisfy the inf-sup condition \Eq{biginfsup-h}. For $[\u_h,p_h] \in V_{h} \times Q_{h}$ taking arbitrary values on $\Gamma$, we have the following result:

\begin{theorem} \label{th-0} 
Suppose that the FE space $V_{h,0}\times Q_{h,0}$ satisfies the inf-sup condition \Eq{biginfsup-h}. Then, \rcr if $L_0 > h$  \ot for each $[\u_h,p_h] \in V_{h} \times Q_{h}$ there exists $[\vt_{h,0},q_{h,0}] \in V_{h,0} \times Q_{h,0}$ such that
\begin{align}
B([\u_h,p_h],[\vt_{h,0} , q_{h,0}]) \gtrsim \Vert [\u_h,p_h] \Vert^2_{V\times Q}   
- \gamma \frac{\nu}{h}\Vert \n\times \u_h\Vert^2_{L^2(\Gamma)}
- \gamma \frac{L_0^2}{\nu h}\Vert p_h\Vert^2_{L^2(\Gamma)},\label{eq:gen-inf-sup}
\end{align}
for a constant $\gamma \geq 0$.
\end{theorem}

\begin{proof} 
Let ${\mathcal T}_\Gamma = \{ K\in {\mathcal T}_h ~\vert ~\partial K \cap\Gamma \not = \emptyset\}$ and  $\Omega_\Gamma = {\rm int}( \bigcup_{K\in {\mathcal T}_\Gamma} K)$, i.e., $\Omega_\Gamma $ is the first layer of element subdomains inside $\Omega$, 
\rc 
and let ${\mathcal T}^0_\Gamma \subset {\mathcal T}_\Gamma $ the subset of elements $K$ such that $\partial K\cap\Gamma$ is an edge if $d=2$ or a face if $d = 3$.
\ot
Let us also write $\partial \Omega_\Gamma  = \Gamma \cup \Gamma_0$. 

Let us consider the splitting $Q_h = Q_{h,0} \oplus Q_{h,\Gamma}$, where $Q_{h,0}$ is the subspace of functions in $Q_h$ vanishing on $\Gamma$ and $Q_{h,\Gamma}$ its complement, i.e., the space made of functions in $Q_h$ which are zero at all the interior nodes of $\Omega$. We may split all functions $q_h\in Q_h$ as $q_h = q_{h,0} + q_{h,\Gamma}$, with $q_{h,0}\in Q_{h,0}$ and $q_{h,\Gamma} \in Q_{h,\Gamma}$. We construct $q_{h,\Gamma} \in Q_{h,\Gamma}$ from the degrees of freedom of $q_{h} \in Q_{h}$ on $\Gamma$ 
and setting to zero all internal degrees of freedom. In particular, $q_{h,\Gamma} = 0$ on $\Gamma_0$.

If $p_h = p_{h,0} + p_{h,\Gamma}$, for $p_{h,\Gamma}$ we have that:
\begin{align}
\Vert \nabla p_{h,\Gamma} \Vert^2_{L^2(\Omega)} & = 
\sum_{K\in {\mathcal{T}_\Gamma}} \Vert \nabla p_{h,\Gamma} \Vert^2_{L^2(K)}
\lesssim \sum_{K\in {\mathcal{T}_\Gamma}} \frac{1}{h^2}  \Vert p_{h,\Gamma} \Vert^2_{L^\infty(K)} h^d \cl
& = \sum_{K\in {\mathcal{T}_\Gamma}} \frac{1}{h^2}  \Vert p_{h,\Gamma} \Vert^2_{L^\infty(\partial K \cap \Gamma)} h^d \cl
&\rc  \lesssim \sum_{K\in {\mathcal{T}^0_\Gamma}} \frac{1}{h^2}  \Vert p_{h,\Gamma} \Vert^2_{L^\infty(\partial K \cap \Gamma)} h^d \ot  \cl
& \lesssim \sum_{K\in {\mathcal{T}^0_\Gamma}} \frac{1}{h^2}  \Vert p_{h} \Vert^2_{L^2(\partial K \cap \Gamma)} h^{-(d-1)} h^d \cl
& \lesssim \frac{1}{h}\Vert p_h\Vert^2_{L^2(\Gamma)}. \label{eq:th3-1-a}
\end{align}
\rc In the fourth step we have used that for elements $K$ such that $\partial K\cap \Gamma$ is a point (or an edge if $d=3$), the norm $\Vert p_{h,\Gamma} \Vert^2_{L^\infty(\partial K \cap \Gamma)}$ is bounded by that of the neighbors that have a whole edge (face, if $d=3$) on $\Gamma$, so that this norm can be absorved by that of the neighbors (and there are a finite number of these elements with the same neighbors if the mesh is non-degenerate). In the fifth step we have used \Eq{invineq-b}. \ot

For $V_h$ we may proceed similarly. Let $V_h = V_{h,0} \oplus V_{h,\Gamma}$, so that each $\vt_h\in V_h$ may be written as $\vt_h = \vt_{h,0} + \vt_{h,\Gamma}$ and $\vt_{h,\Gamma} \in V_{h,\Gamma}$ is constructed such that $\n\times \vt_{h,\Gamma} = \n\times \vt_{h}$ on $\Gamma$, i.e., the degrees of freedom associated to $\Gamma$ of $\n\times \vt_{h,\Gamma}$ are set equal to those of $\n\times \vt_{h}$, 
and all internal degrees of freedom of $\vt_{h,\Gamma}$ are zero. If needed, we also set $\vt_{h,\Gamma}\cdot\n = 0$ and, if $d=3$, $\n\times\vt_{h,\Gamma}\times\n = {\bf 0}$ on $\Gamma$. 

If $\u_h = \u_{h,0} + \u_{h,\Gamma}$, for $\u_{h,\Gamma}$ we have that:
\begin{align}
\Vert \nabla \times \u_{h,\Gamma} \Vert^2_{L^2(\Omega)} & = 
\sum_{K\in {\mathcal{T}_\Gamma}} \Vert \nabla \times \u_{h,\Gamma} \Vert^2_{L^2(K)}
\lesssim \sum_{K\in {\mathcal{T}_\Gamma}} \frac{1}{h^2}  \Vert \u_{h,\Gamma} \Vert^2_{L^\infty(K)} h^d \cl
& = \sum_{K\in {\mathcal{T}_\Gamma}} \frac{1}{h^2}  \Vert \u_{h,\Gamma} \Vert^2_{L^\infty(\partial K \cap \Gamma)} h^d 
\rc  \lesssim \sum_{K\in {\mathcal{T}^0_\Gamma}} \frac{1}{h^2}  \Vert
\n\times\ot\db \u_{h,\Gamma}\ot\rc \Vert^2_{L^\infty(\partial K \cap \Gamma)} h^d \ot \cl
& \lesssim \sum_{K\in {\mathcal{T}^0_\Gamma}} \frac{1}{h^2}  \Vert \n\times\u_{h} \Vert^2_{L^2(\partial K \cap \Gamma)} h^{-(d-1)} h^d \cl
& \lesssim \frac{1}{h}\Vert \n\times\u_h\Vert^2_{L^2(\Gamma)}.\label{eq:th3-1-b}
\end{align}
Using similar arguments we easily get that
\begin{align}
\Vert \u_{h,\Gamma} \Vert^2_{L^2(\Omega)} \lesssim {h}\Vert \n\times\u_h\Vert^2_{L^2(\Gamma)}.\label{eq:th3-1-c}
\end{align}

Given $[\u_{h,0}, p_{h,0}] = [\u_{h}, p_{h}]  - [\u_{h,\Gamma}, p_{h,\Gamma}]
\in V_{h,0}\times Q_{h,0}$ constructed as explained above, let $[\vt_{h,0},
q_{h,0}]\in V_{h,0}\times Q_{h,0}$ be the element for which \Eq{biginfsup-h} holds. We then have:
\begin{align}
B([\u_h,p_h],[\vt_{h,0} , q_{h,0}]) & \geq K_B \Vert [\u_{h,0}, p_{h,0}] \Vert_{V\times Q}\Vert [\vt_{h,0}, q_{h,0}] \Vert_{V\times Q} \cl
& + \nu (\nabla\times \u_{h,\Gamma} , \nabla\times \vt_{h,0})_\Omega + (\vt_{h,0} , \nabla p_{h,\Gamma} )_\Omega + (\u_{h,\Gamma} , \nabla q_{h,0} )_\Omega.\label{eq:th3-1-cc}
\end{align}
In the following, $\alpha_i > 0$, $i= 1,2,3$, denote constants arising from Young's inequality. Choosing $\Vert [\vt_{h,0}, q_{h,0}] \Vert_{V\times Q} = \Vert [\u_{h,0}, p_{h,0}] \Vert_{V\times Q}$, we obtain:
\begin{align*}
B([\u_h,p_h],[\vt_{h,0} & , q_{h,0}])  \geq K_B \Vert [\u_{h,0}, p_{h,0}] \Vert^2_{V\times Q}\cl
& - \frac{\nu}{2\alpha_1} \Vert \nabla\times \u_{h,\Gamma} \Vert^2_{L^2(\Omega)} 
- \frac{L_0^2}{2\alpha_2 \nu} \Vert \nabla p_{h,\Gamma} \Vert^2_{L^2(\Omega)} 
- \frac{\nu}{2\alpha_3 L_0^2}  \Vert\u_{h,\Gamma} \Vert^2_{L^2(\Omega)}  \cl
& - \frac{1}{2}(\alpha_1 + \alpha_2 + \alpha_3)  \Vert [\u_{h,0}, p_{h,0}] \Vert^2_{V\times Q}
\end{align*}
Taking $\alpha_i$ sufficiently small, $i=1,2,3$, and making use of \Eq{th3-1-a}-\Eq{th3-1-c} we get, assuming $L_0 > h = \max_K \{h_K\}$:
\begin{align}
B([\u_h,p_h],[\vt_{h,0} , q_{h,0}]) \gtrsim \Vert [\u_{h,0},p_{h,0}] \Vert^2_{V\times Q}   
- \gamma_0 \frac{\nu}{h}\Vert \n\times \u_h\Vert^2_{L^2(\Gamma)}
- \gamma_0 \frac{L_0^2}{\nu h}\Vert p_h\Vert^2_{L^2(\Gamma)},\label{eq:th3-1-d}
\end{align}
for a constant $\gamma_0 \geq 0$. Furthermore, using again \Eq{th3-1-a}-\Eq{th3-1-c} we get:
\begin{align*}
\Vert [\u_{h,0},p_{h,0}] \Vert^2_{V\times Q}   
& \gtrsim \Vert [\u_{h},p_{h}] \Vert^2_{V\times Q}  - \Vert [\u_{h,\Gamma},p_{h,\Gamma}] \Vert^2_{V\times Q}   \cl
& \gtrsim \Vert [\u_{h},p_{h}] \Vert^2_{V\times Q}  
- \gamma_1 \frac{\nu}{h}\Vert \n\times \u_h\Vert^2_{L^2(\Gamma)}
- \gamma_1 \frac{L_0^2}{\nu h}\Vert p_h\Vert^2_{L^2(\Gamma)},
\end{align*}
for a constant $\gamma_1 \geq 0$, which combined with \Eq{th3-1-d} yields the theorem.
\end{proof}

Estimate \Eq{gen-inf-sup} explicitly displays which terms spoil stability of the problem without boundary conditions. The terms introduced by Nitsche's method need precisely to compensate them.

If no boundary conditions are prescribed, from identity \Eq{byparts} it is found that the discrete weak form of the differential equation \Eq{m1} would be
\begin{align}
B([\u_h,p_h] , [\vt_h , q_h]) - \langle {\mathcal F}_n ( [\u_h , p_h]) , {\mathcal D} ([\vt_h , q_h])\rangle_\Gamma =  \langle \vt_h , {\bm f}\rangle_\Omega,\label{eq:nit1}
\end{align}
the different terms being defined in \Eq{defB}-\Eq{defD}. For the continuous solution $[\u,p]\in V\times Q$, there holds
\begin{align*}
- \langle {\mathcal F}_n ( [\vt_h , q_h]) , {\mathcal D} ([\u , p])\rangle_\Gamma & = - \langle {\mathcal F}_n ( [\vt_h , q_h]) , {\mathcal D} ([\bar{\u} , 0])\rangle_\Gamma,\\
  \langle {\mathcal D} ( [\vt_h , q_h]) , {\bm N}{\mathcal D} ([\u, p])\rangle_\Gamma & =  \langle {\mathcal D} ( [\vt_h , q_h]) , {\bm N}{\mathcal D} ([\bar{\u} , 0])\rangle_\Gamma,
\end{align*}
where ${\bm N}$ is a matrix that scales the vector of Dirichlet boundary conditions. The symmetric version of Nitsche's method we use is obtained by adding to Eq. \Eq{nit1} these two expressions evaluated with the FE solution $[\u_h,p_h]\in V_h\times Q_h$. Taking the scaling matrix as
\begin{align*}
{\bm N} = {\rm diag} \left(N_u \frac{\nu}{h} {\bm I} , N_p \frac{L_0^2}{h}\right), 
\end{align*}
where $N_u$ and $N_p$ are dimensionless algorithmic constants that need to be determined, the final problem is: find $[\u_h,p_h]\in V_h\times Q_h$ such that 
\begin{align}
B_{\rm N} ([\u_h,p_h] , [\vt_h,q_h]) = L_{\rm N}  ([\vt_h,q_h])\qquad \forall [\vt_h,q_h] \in V_h\times Q_h, \label{eq:gal-nit}
\end{align}
where
\begin{align}
B_{\rm N} ([\u_h,p_h] , [\vt_h,q_h]) 
& = \nu ( \nabla\times \vt_h , \nabla\times \u_h)_\Omega + (\vt_h ,  \nabla p_h)_\Omega + ( \u_h , \nabla q_h)_\Omega \cl
& \quad- \nu \langle {\bm n} \times \vt_h , \nabla\times \u_h\rangle_\Gamma - \langle{\bm  n} \cdot \u_h , q_h\rangle_\Gamma \cl
& \quad  - \nu \langle{\bm n}  \times \u_h , \nabla\times \vt_h\rangle_\Gamma - \langle {\bm n}  \cdot \vt_h , p_h\rangle_\Gamma \cl
& \quad + N_u \frac{\nu}{h} \langle {\bm n}  \times \vt_h , {\bm n}  \times \u_h \rangle_\Gamma 
- N_p \frac{L_0^2}{\nu h} ( p_h , q_h)_\Gamma  \label{eq:gal-nit-B}\\
L_{\rm N}  ([\vt_h,q_h])
& = \langle \vt_h , {\bm f}\rangle_\Omega- \nu \langle{\bm n}  \times \bar{\u} , \nabla\times \vt_h\rangle_\Gamma 
+ N_u \frac{\nu}{h} \langle {\bm n}  \times \vt_h , {\bm n} \times \bar{\u} \rangle_\Gamma. \label{eq:gal-nit-L}
\end{align}

Apart from the boundary term arising from integration by parts, the rest of boundary terms introduced (symmetrisation and penalisation) can be interpreted as stabilisation terms provided by the sub-grid scales on the boundary in the context of the variational multi-scale method. This interpretation is introduced in \cite{codina-et-al-2009-1}.

\begin{theorem}\label{th:gal-nit} 
Assume that the FE space $V_{h,0}\times Q_{h,0}$ satisfies the inf-sup condition \Eq{biginfsup-h}. Then, for $N_u$ and $N_p$ sufficiently large, $B_{\rm N}$ is inf-sup stable in the norm
\begin{align}
\vtres{[\vt_h, q_h]}^2_{V\times Q,{\rm N}} = \Vert [\vt_h, q_h] \Vert^2_{V\times Q}  
+ \frac{\nu}{h} \Vert \n \times \vt_h \Vert^2_{L^2(\Gamma)}
+ \frac{L_0^2}{\nu h} \Vert q_h  \Vert^2_{L^2(\Gamma)}. \el
\end{align}
\end{theorem}

\begin{proof}
It is trivially checked that $\vtres{\cdot}_{V\times Q,{\rm N}}$ is indeed a norm in $V_h\times Q_h$.

Let us start noting that
\begin{align*}
B_{\rm N} ([\u_h,p_h] , [\u_h,-p_h]) 
& \geq \nu \Vert \nabla \times \u_h \Vert^2_{L^2(\Omega)}
- 2 \nu \Vert \n \times \u_h \Vert_{L^2(\Gamma)} \Vert \nabla \times \u_h \Vert_{L^2(\Gamma)} \cl
& + N_u \frac{\nu}{h} \Vert \n \times \u_h \Vert^2_{L^2(\Gamma)} 
+ N_p \frac{L_0^2}{\nu h} \Vert p_h  \Vert^2_{L^2(\Gamma)}.
\end{align*}
Using the trace inequality \Eq{traceineq} and Young's inequality we get, for all $\alpha > 0$:
\begin{align*}
B_{\rm N} ([\u_h,p_h] , [\u_h,-p_h]) 
& \geq \nu \Vert \nabla \times \u_h \Vert^2_{L^2(\Omega)} 
+ N_u \frac{\nu}{h} \Vert \n \times \u_h \Vert^2_{L^2(\Gamma)}
+ N_p \frac{L_0^2}{\nu h} \Vert p_h  \Vert^2_{L^2(\Gamma)} \cl
& - 2\nu \left( \frac{1}{2\alpha} \Vert \n \times \u_h \Vert^2_{L^2(\Gamma)} + \frac{\alpha C^2_{\rm trace}}{2h} \Vert \nabla \times \u_h \Vert^2_{L^2(\Omega)} \right).\el
\end{align*}
Taking for example $\alpha = h (2C^2_{\rm trace})^{-1}$ and assuming $N_u \geq 2C^2_{\rm trace} + N_u'$, with $N_u' > 0$:
\begin{align}
B_{\rm N} ([\u_h,p_h] , [\u_h,-p_h]) 
& \geq \frac{\nu}{2} \Vert \nabla \times \u_h \Vert^2_{L^2(\Omega)} 
+  N_u' \frac{\nu}{h} \Vert \n \times \u_h \Vert^2_{L^2(\Gamma)}
+ N_p \frac{L_0^2}{\nu h} \Vert p_h  \Vert^2_{L^2(\Gamma)}.\label{eq:th1-1}
\end{align}

Let now $[\vt_{h,0},q_{h,0}]\in V_{h,0}\times Q_{h,0}$ be the pair whose existence is established in Theorem~\ref{th-0} that satisfies \Eq{gen-inf-sup}. Recall that $\Vert [\vt_{h,0},q_{h,0}]\Vert_{V\times Q} = \Vert [\u_{h,0},q_{h,0}]\Vert_{V\times Q} \leq \Vert [\u_{h},p_{h}]\Vert_{V\times Q}$. Since $\n\times \vt_{h,0} = {\bf 0}$ and $q_{h,0} = 0$ on $\Gamma$, we have that
\begin{align}
B_{\rm N}([\u_h,p_h]  & ,[\vt_{h,0} , q_{h,0}])
  = B([\u_h,p_h],[\vt_{h,0} , q_{h,0}])\cl
&  - \nu \langle{\bm n}  \times \u_h , \nabla\times \vt_{h,0}\rangle_\Gamma - \langle {\bm n}  \cdot \vt_{h,0} , p_h\rangle_\Gamma\cl
& \gtrsim \Vert [\u_h,p_h] \Vert^2_{V\times Q}   
- \gamma \frac{\nu}{h}\Vert \n\times \u_h\Vert^2_{L^2(\Gamma)}
- \gamma \frac{L_0^2}{\nu h}\Vert p_h\Vert^2_{L^2(\Gamma)} \cl
& - \nu \langle{\bm n}  \times \u_h , \nabla\times \vt_{h,0}\rangle_\Gamma - \langle {\bm n}  \cdot \vt_{h,0} , p_h\rangle_\Gamma. \label{eq:th3-2-a}
\end{align}
We may now bound the last two terms as follows:
\begin{align}
 \nu \langle{\bm n}  \times \u_h , \nabla\times \vt_{h,0}\rangle_\Gamma 
&  \leq \nu \Vert \nabla \times \vt_{h,0} \Vert_{L^2(\Gamma)}\Vert \n \times \u_h \Vert_{L^2(\Gamma)}  \cl
& \leq \nu \frac{C_{\rm trace}}{h^{1/2}} \Vert \nabla \times \vt_{h,0} \Vert_{L^2(\Omega)}\Vert \n \times \u_h \Vert_{L^2(\Gamma)} \cl
& \leq \frac{\alpha_1}{2} \Vert [\u_{h},p_{h}]\Vert^2_{V\times Q}  + \frac{1}{2\alpha_1} {C^2_{\rm trace}} \frac{\nu}{h} \Vert \n \times \u_h \Vert^2_{L^2(\Gamma)}, \el
\end{align}
and 
\begin{align}
\langle {\bm n}  \cdot \vt_{h,0} , p_h\rangle_\Gamma 
& \leq \Vert \n\cdot \vt_{h,0} \Vert_{L^2(\Gamma)}\Vert p_h \Vert_{L^2(\Gamma)} \cl
& \leq \frac{C_{\rm trace}}{h^{1/2}}  \Vert  \vt_{h,0} \Vert_{L^2(\Omega)} \Vert p_h \Vert_{L^2(\Gamma)} \cl
& \leq \frac{\alpha_2}{2} \Vert [\u_{h},p_{h}]\Vert^2_{V\times Q}  + \frac{1}{2\alpha_2} {C^2_{\rm trace}} \frac{L_0^2}{\nu h} \Vert p_h \Vert^2_{L^2(\Gamma)}.\el
\end{align}
For $\alpha_1$ and $\alpha_2$ small enough, \Eq{th3-2-a} yields
\begin{align}
B_{\rm N}([\u_h,p_h],[\vt_{h,0} , q_{h,0}])
& \gtrsim \Vert [\u_h,p_h] \Vert^2_{V\times Q}   
- \gamma^\ast \frac{\nu}{h}\Vert \n\times \u_h\Vert^2_{L^2(\Gamma)}
- \gamma^\ast\frac{L_0^2}{\nu h}\Vert p_h\Vert^2_{L^2(\Gamma)}, \label{eq:th3-2-b}
\end{align}
for a constant $\gamma^\ast > 0$. 

Set now $[\vt_h,q_h] = [\u_h , -p_h] + \delta[\vt_{h,0} , q_{h,0}]$, with $\delta$ small enough (or $N'_u$ and $N_p$ large enough). Combining \Eq{th1-1} and  \Eq{th3-2-b} it follows that
\begin{align}
B_{\rm N}([\u_h,p_h],[\vt_{h} , q_{h}]) \gtrsim \vtres{[\u_h, p_h]}^2_{V\times Q,{\rm N}}.
\end{align}
The proof concludes after checking that  $ \vtres{[\vt_h, q_h]}_{{V\times Q},{\rm N}} \leq (1+\delta) \vtres{[\u_h, p_h]}_{V\times Q,{\rm N}}$.
\end{proof}

Let us prove two preliminary results to obtain the analogous of Theorem~\ref{th:th1} when Dirichlet conditions are prescribed using Nitsche's method:

\begin{lemma}\label{lem:l1}
The linear form $L_{\rm N}$ given in \Eq{gal-nit-L} is continuous in the norm $\vtres{\cdot}_{V\times Q,{\rm N}}$, the continuity constant being bounded as
\begin{align}
\vtres{L_{\rm N}}_{{\mathcal L}(V_h\times Q_h,{\rm N}; \mathbb{R})} \lesssim  \Vert {\bm f}\Vert_{V'} +  \Bigl(\frac{\nu}{h}\Bigr)^{1/2} \Vert \n\times \bar{\u}\Vert_{L^2(\Gamma)}.\label{eq:ineq-l31}
\end{align}
\end{lemma}

\begin{proof}
For any $[\vt_h,q_h]\in V_h\times Q_h$ we have that:
\begin{align*}
L_{\rm N}([\vt_h,q_h] ) 
&\lesssim \Vert {\bm f} \Vert_{V'} \Vert \vt_h \Vert_V 
+ \nu \Vert \n \times \bar{\u} \Vert_{L^2(\Gamma)} \Vert \nabla \times \vt_h \Vert_{L^2(\Gamma)} \cl
& + \frac{\nu}{h}\Vert \n \times \bar{\u} \Vert_{L^2(\Gamma)} \Vert \n \times \vt_h \Vert_{L^2(\Gamma)} \cl
& \lesssim \Vert {\bm f} \Vert_{V'} \Vert \vt_h \Vert_V 
+ \frac{\nu^{1/2}}{h^{1/2}}  \Vert \n \times \bar{\u} \Vert_{L^2(\Gamma)} C_{\rm trace} \nu^{1/2} \Vert \nabla \times \vt_h \Vert_{L^2(\Omega)}\cl
& + \frac{\nu^{1/2}}{h^{1/2}}  \Vert \n \times \bar{\u} \Vert_{L^2(\Gamma)} \frac{\nu^{1/2}}{h^{1/2}}  \Vert \n \times {\vt_h} \Vert_{L^2(\Gamma)} \cl
& \lesssim 
\Bigl( \Vert {\bm f}\Vert_{V'} + \frac{\nu^{1/2}}{h^{1/2}} \Vert \n\times \bar{\u}\Vert_{L^2(\Gamma)}\Bigr) \vtres{[ \vt_h, q_h] }_{V\times Q,{\rm N}},
\end{align*}
thus proving the Lemma.
\end{proof}

\begin{lemma}\label{lem:l2}  
For any $[\u,p]\in V\times Q$, let the interpolation error function be
\begin{align}
E(\u,p;h) = \inf_{[ \tilde{\u}_h, \tilde{p}_h] \in V_{h}\times Q_{h}} D([ \u - \tilde{\u}_h, p - \tilde{p}_h]),\label{eq:err-est-gal-nit}
\end{align}
where
\begin{align*}
D([ \vt ,  q])
 = \vtres{ [\vt , q ] }_{V\times Q,{\rm N}}  
+ {(\nu h)^{1/2}} \Vert \n\times \nabla \times  \vt \Vert_{L^2(\Gamma)} 
+  \frac{(\nu h)^{1/2}}{L_0} \Vert \n\cdot \vt \Vert_{L^2(\Gamma)}.
\end{align*}
Then, for all $[\vt_h,q_h]\in V_h\times Q_h$ there holds
\begin{align}
\inf_{[ \tilde{\u}_h, \tilde{p}_h] \in V_{h}\times Q_{h}}  B_{\rm N}( [\u - \tilde{\u}_h, p - \tilde{p}_h] , [\vt_h,q_h]) \lesssim E(\u,p;h) \vtres{[\vt_h,q_h]}_{V\times Q,{\rm N}}. \label{eq:ineq-l32}
\end{align}
\end{lemma}

\begin{proof} The terms involving volume integrals and the penalisation terms in $B_{\rm N}( [\u - \tilde{\u}_h, p - \tilde{p}_h] , [\vt_h,q_h]) $ are bounded by $  \vtres{[\u - \tilde{\u}_h, p - \tilde{p}_h]}_{V\times Q,{\rm N}}\vtres{[\vt_h,q_h]}_{V\times Q,{\rm N}}$, as it is immediately checked. For the rest of boundary terms we can proceed as follows:
\begin{align*}
& - \nu \langle {\bm n} \times \vt_h , \nabla\times (\u - \tilde{\u}_h)\rangle_\Gamma   \cl
& \qquad \lesssim \Bigl(\frac{\nu}{h}\Bigr)^{1/2} \Vert  {\bm n} \times \vt_h \Vert_{L^2(\Gamma)} {(\nu h)^{1/2}}  \Vert  {\bm n} \times \nabla \times (\u - \tilde{\u}_h) \Vert_{L^2(\Gamma)},
\cl
& - \langle{\bm  n} \cdot (\u - \tilde{\u}_h) , q_h\rangle_\Gamma
\lesssim \frac{L_0}{(\nu h)^{1/2}} \Vert q_h \Vert_{L^2(\Gamma)}  \frac{(\nu h)^{1/2}}{L_0} \Vert \n\cdot (\u - \tilde{\u}_h) \Vert_{L^2(\Gamma)},
\cl
& - \nu \langle{\bm n}  \times (\u - \tilde{\u}_h) , \nabla\times \vt_h\rangle_\Gamma \cl
& \qquad \lesssim C_{\rm trace}  \nu^{1/2} \Vert  \nabla\times \vt_h\Vert_{L^2(\Omega)}\Bigl(\frac{\nu}{h}\Bigr)^{1/2} \Vert \n\times  (\u - \tilde{\u}_h) \Vert_{L^2(\Gamma)},
\cl
& - \langle {\bm n}  \cdot \vt_h , p - \tilde{p}_h\rangle_\Gamma
\lesssim C_{\rm trace} \frac{{\nu}^{1/2}}{L_0} \Vert \vt_h\Vert_{L^2(\Omega)} \frac{L_0}{(\nu h)^{1/2}} \Vert p - \tilde{p}_h \Vert_{L^2(\Gamma)}.
\el
\end{align*}
Clearly, all these terms are bounded by $D([\u -\tilde{\u}_h, p - \tilde{p}_h]) \vtres{[\vt_h,q_h]}_{V\times Q,{\rm N}}$, from where
\begin{align}
B_{\rm N}( [\u - \tilde{\u}_h, p - \tilde{p}_h] , [\vt_h,q_h]) \leq D([ \u - \tilde{\u}_h, p - \tilde{p}_h]) \vtres{[\vt_h,q_h]}_{V\times Q,{\rm N}},\label{eq:l2tem}
\end{align}
and the result follows taking the infimum for $[\tilde{\u}_h, \tilde{p}_h]$ over $V_h\times Q_h$.
\end{proof}

Let $k_u$ be the highest order of the complete piecewise polynomial contained in $V_h$ and $k_p$ the one of the complete piecewise polynomial contained in $Q_h$. Using standard interpolation estimates, it is seen that
\begin{align}
E(\u,p;h) \lesssim \nu^{1/2} h^{s_u-1} \Vert \u \Vert_{H^{s_u}(\Omega)} + \frac{L_0}{\nu^{1/2}} h^{s_p-1} \Vert p \Vert_{H^{s_p}(\Omega)},\label{eq:asbeh}
\end{align}
where \os $s_u = \min\{ k_u+1, r_u\}$, $s_p = \min\{ k_p + 1, r_p\}$ \ot and $r_u$ and $r_p$ are the Sobolev regularity of $\u\in V$ and $p \in Q$, respectively. If we prove that this is the error function of the formulation, it will be clearly optimal. This is indeed proved in the following result, which is the analogous of Theorem~\ref{th:th1} when Dirichlet conditions are prescribed using Nitsche's method:

\begin{theorem}\label{th:conv-gal-nit} 
Under the assumptions of Theorem~\ref{th:gal-nit}, problem \Eq{gal-nit} is well posed, in the sense that it admits a unique solution $[\u_h,p_h]\in V_{h}\times Q_{h}$ that satisfies
\begin{align}
\vtres{ [\u_h,p_h] \Vert}_{V\times Q,{\rm N}} \lesssim \Vert {\bm f}\Vert_{V'}
+ \Bigl(\frac{\nu}{h}\Bigr)^{1/2} \Vert \n\times
\bar{\u}\Vert_{L^2(\Gamma)}.\label{eq:gal-nit-unstab}
\end{align}
Furthermore, $[\u_h,p_h]$ converges optimally as $h \to 0$ to the solution $[\u,p]\in V \times Q$ of the continuous problem \Eq{varc1}-\Eq{varc2}, in the following sense:
\begin{align}
& \vtres{ [\u - \u_h, p - p_h] }_{V\times Q,{\rm N}} \lesssim  E(\u,p; h),\el
\end{align}
where $E(\u,p; h)$ is given in \Eq{err-est-gal-nit}.
\end{theorem}

\begin{proof}
Existence and uniqueness of the discrete solutions follows from the inf-sup
condition stated in Theorem~\ref{th:gal-nit}, and the stability estimate
\Eq{gal-nit-unstab} is a direct consequence of the continuity of $L_{\rm N}$ proved in Lemma~\ref{lem:l1} and the inf-sup condition. 

As discussed above, problem \Eq{gal-nit} is consistent, that is to say,  $B_{\rm N} ([\u,p] , [\vt_h,q_h]) = L_{\rm N}  ([\vt_h,q_h])$ for all $[\vt_h,q_h] \in V_h\times Q_h$, and therefore $B_{\rm N} ([\u - \u_h,p - p_h] , [\vt_h,q_h]) = 0$ for all $[\vt_h,q_h] \in V_h\times Q_h$. Let us pick $[\tilde{\u}_h, \tilde{p}_h] \in V_h\times Q_h$, arbitrary. Now the proof is standard:
\begin{alignat*}{3}
& \vtres{[\u_h - \tilde{\u}_h,p_h -  \tilde{p}_h ] }_{V\times Q,{\rm N}} \vtres{[\vt_h,q_h]}_{V\times Q,{\rm N}}  \cl
& \qquad \lesssim B([\u_h -\tilde{\u}_h, p_h - \tilde{p}_h] , [\vt_h,q_h])  &&\quad \hbox{from the inf-sup condition,}\cl
& \qquad = B([\u -\tilde{\u}_h, p - \tilde{p}_h] , [\vt_h,q_h])  &&\quad  \hbox{from consistency,} \cl
& \qquad \lesssim D([\u -\tilde{\u}_h, p - \tilde{p}_h]) \vtres{[\vt_h,q_h]}_{V\times Q,{\rm N}}  && \quad \hbox{from \Eq{l2tem}.}
\end{alignat*}
Since $\vtres{[\u - \tilde{\u}_h,p -  \tilde{p}_h ] }_{V\times Q,{\rm N}} \leq D([\u -\tilde{\u}_h, p - \tilde{p}_h])$, the theorem follows from the triangle inequality and taking the infimum for $[\tilde{\u}_h, \tilde{p}_h]$ over $V_h\times Q_h$. 
\end{proof}

\rc
Clearly, the constant involved in inequality \eqref{eq:th3-2-b} is independent of $N_u$ and $N_p$ when they are large; more precisely, this constant behaves as $\min \{ 1, N_u, N_p\}$. Thus, the inf-sup constant in the inf-sup condition stated in Theorem~\ref{th:gal-nit} is bounded as $N_u,N_p\to\infty$. On the contrary, the constants involved in inequalities \eqref{eq:ineq-l31} and \eqref{eq:ineq-l32} grow as $\max\{ 1 , N_u , N_p\}$ when $N_u,N_p\to\infty$. As a consequence, the error estimate provided by Theorem~\ref{th:conv-gal-nit} grows as $\max\{ N_u , N_p\}$. In practice it is convenient to take these algorithmic constants as small as possible, although large enough to fulfil the requirements found in the proof of Theorem~\ref{th:gal-nit}.
\ot

\subsection{Nitsche's method using the stabilised FE approximation}
\label{subsec:nitschestabilised}
We consider now Nitsche's method in combination with the stabilised formulation presented in section~\ref{sec:stab-exactBC}. The analysis is similar to that of the Galerkin method, and therefore we will only concentrate on the minor differences introduced by the stabilising terms.

Let us start with the counterpart of Theorem~\ref{th-0}: 

\begin{theorem} \label{th-0-stab} 
Consider the stabilised bilinear form \Eq{defBS}. Then, for each $[\u_h,p_h] \in V_{h} \times Q_{h}$ there exists $[\vt_{h,0},q_{h,0}] \in V_{h,0} \times Q_{h,0}$ such that
\begin{align}
B_{\rm S}([\u_h,p_h],[\vt_{h,0} , q_{h,0}]) \gtrsim \Vert [\u_h,p_h] \Vert^2_{V\times Q,{\rm S}}   
- \gamma \frac{\nu}{h}\Vert \n\times \u_h\Vert^2_{L^2(\Gamma)}
- \gamma \frac{L_0^2}{\nu h}\Vert p_h\Vert^2_{L^2(\Gamma)},\label{eq:gen-inf-sup-stab}
\end{align}
for a constant $\gamma \geq 0$.
\end{theorem}

\begin{proof}
The proof is very similar to that of Theorem~\ref{th-0}. In particular, given $\u_h$, $\u_{h,\Gamma}$ is constructed in the same way as in Theorem~\ref{th-0}, as we wish that $\u_{h,0} = \u_h - \u_{h,\Gamma}$ satisfies that $\n\times \u_{h,0} = {\bf 0}$ on $\Gamma$. 

By virtue of Theorem~\ref{th:inf-sup-BS}, now we will obtain, instead of \Eq{th3-1-c}:
\begin{align}
B_{\rm S}([\u_h,p_h],[\vt_{h,0} , q_{h,0}]) & \gtrsim \Vert [\u_{h,0}, p_{h,0}] \Vert_{V\times Q,{\rm S}}\Vert [\vt_{h,0}, q_{h,0}] \Vert_{V\times Q,{\rm S}} \cl
& + \nu (\nabla\times \u_{h,\Gamma} , \nabla\times \vt_{h,0})_\Omega + (\vt_{h,0} , \nabla p_{h,\Gamma} )_\Omega + (\u_{h,\Gamma} , \nabla q_{h,0} )_\Omega\cl
& + c_u \frac{h^2 \nu}{L_0^2} (\nabla \cdot \u_{h,\Gamma} , \nabla \cdot \vt_{h,0})_\Omega - \frac{L_0^2}{\nu} ( \nabla p_{h,\Gamma} , \nabla p_{h,0})_\Omega,\el
\end{align}
for a certain $[\vt_{h,0} , q_{h,0}]\in V_{h,0}\times Q_{h,0}$. The meaning of different variables and unknowns is the same as in Theorem~\ref{th:inf-sup-BS}. Now we have to deal with the last two terms of this expression, which offer no difficulty, as: 
\begin{align}
& c_u \frac{h^2 \nu}{L_0^2} (\nabla \cdot \u_{h,\Gamma} , \nabla \cdot \vt_{h,0})_\Omega 
\lesssim  \frac{1}{2}\alpha_4 \Vert [\u_{h,0} , p_{h,0}]\Vert_{V\times Q,{\rm S}} + \frac{\nu h^2 }{2\alpha_4 L_0^2} \Vert \nabla\cdot\u_{h,\Gamma}\Vert_{L^2(\Omega)},\cl
& \frac{L_0^2}{\nu} ( \nabla p_{h,\Gamma} , \nabla p_{h,0})_\Omega 
 \lesssim  \frac{1}{2}\alpha_5 \Vert [\u_{h,0} , p_{h,0}]\Vert_{V\times Q,{\rm S}} + \frac{L_0^2}{2\alpha_5\nu}\Vert  \nabla p_{h,\Gamma}\Vert_{L^2(\Omega)}.\el
\end{align}
Using the same steps as in \Eq{th3-1-b} it is easily checked that 
\begin{align}
\Vert \nabla \cdot \u_{h,\Gamma} \Vert^2_{L^2(\Omega)}  \lesssim \frac{1}{h}\Vert \n\times\u_h\Vert^2_{L^2(\Gamma)},\el
\end{align}
and we already proved \Eq{th3-1-a}. The proof concludes as that of Theorem~\ref{th-0}.
\end{proof}

According to this result, the terms that need to be compensated to get stability using Nitsche's method are the same as for the inf-sup stable case. Using the general idea described in section~\ref{sec:gal-nit}, this method reads as follows: find $[\u_h,p_h]\in V_h\times Q_h$ such that 
\begin{align}
B_{\rm SN} ([\u_h,p_h] , [\vt_h,q_h]) = L_{\rm N}  ([\vt_h,q_h])\qquad \forall [\vt_h,q_h] \in V_h\times Q_h, \label{eq:stab-nit}
\end{align}
where
\begin{align}
& B_{\rm SN} ([\u_h,p_h] , [\vt_h,q_h]) 
 = B_{\rm S} ([\u_h,p_h] , [\vt_h,q_h]) \cl
& \qquad- \nu \langle {\bm n} \times \vt_h , \nabla\times \u_h\rangle_\Gamma - \langle{\bm  n} \cdot \u_h , q_h\rangle_\Gamma 
  - \nu \langle{\bm n}  \times \u_h , \nabla\times \vt_h\rangle_\Gamma - \langle {\bm n}  \cdot \vt_h , p_h\rangle_\Gamma \cl
& \qquad + \frac{L_0^2}{\nu }\langle {\bm n} \cdot \nabla p_h , q_h\rangle_\Gamma
+ \frac{L_0^2}{\nu }\langle p_h , \n \cdot \nabla q_h\rangle_\Gamma
  + N_u \frac{\nu}{h} \langle {\bm n}  \times \vt_h , {\bm n}  \times \u_h \rangle_\Gamma \cl
& \qquad - N_p \frac{L_0^2}{\nu h} ( p_h , q_h)_\Gamma  \cl
& \quad = B_{\rm N} ([\u_h,p_h] , [\vt_h,q_h]) 
+ c_u\frac{\nu h^2}{L_0^2} (\nabla\cdot\u_h , \nabla\cdot\vt_h)_\Omega- \frac{L_0^2}{\nu }(\nabla p_h , \nabla q_h)_\Omega\cl
& \qquad+ \frac{L_0^2}{\nu }\langle {\bm n} \cdot \nabla p_h , q_h\rangle_\Gamma
+ \frac{L_0^2}{\nu }\langle p_h , \n \cdot \nabla q_h\rangle_\Gamma.
\label{eq:stab-nit-B}
\end{align}
The first expression corresponds to adding to the stabilised bilinear form Nitsche's terms and the second to adding to the Nitsche's form of the Galerkin method the stabilisation terms and the boundary term arising from the integration by parts of the Laplacian of $p$ and its symmetric counterpart. In fact, since the exact solution is $p=0$, these last two terms could be removed from the formulation.

The analysis proceeds as for the Galerkin case. Let us start with the analogous to Theorem~\ref{th:gal-nit}:

\begin{theorem}\label{th:stab-nit} 
Consider the stabilised bilinear form using Nitsche's method given by \Eq{stab-nit-B}. Then, for $N_u$ and $N_p$ sufficiently large, $B_{\rm SN}$ is inf-sup stable in the norm
\begin{align}
\vtres{[\vt_h, q_h]}^2_{V\times Q,{\rm SN}} = \Vert [\vt_h, q_h] \Vert^2_{V\times Q,{\rm S}}  
+ \frac{\nu}{h} \Vert \n \times \vt_h \Vert^2_{L^2(\Gamma)}
+ \frac{L_0^2}{\nu h} \Vert q_h  \Vert^2_{L^2(\Gamma)}, \el
\end{align}
\end{theorem}

\begin{proof}
One can follow the same steps as in the proof of Theorem~\ref{th:gal-nit}. Again, it is trivially checked that $\vtres{\cdot}_{V\times Q,{\rm SN}}$ is a norm in $V_h\times Q_h$.

Now we have that
\begin{align*}
& B_{\rm SN} ([\u_h,p_h] , [\u_h,-p_h]) 
 \geq \nu \Vert \nabla \times \u_h \Vert^2_{L^2(\Omega)} + \frac{L_0^2}{\nu} \Vert \nabla p_h\Vert^2_{L^2(\Omega)} + \frac{\nu h^2}{L_0^2 } \Vert \nabla\cdot\u_h\Vert^2_{L^2(\Omega)} \cl
& \qquad - 2 \nu \Vert \n \times \u_h \Vert_{L^2(\Gamma)} \Vert \nabla \times \u_h \Vert_{L^2(\Gamma)} 
- 2 \frac{L_0^2}{\nu} \Vert p_h\Vert_{L^2(\Gamma)} \Vert \n\cdot \nabla p_h\Vert_{L^2(\Gamma)}\cl
& \qquad 
+ N_u \frac{\nu}{h} \Vert \n \times \u_h \Vert^2_{L^2(\Gamma)} 
+ N_p \frac{L_0^2}{\nu h} \Vert p_h  \Vert^2_{L^2(\Gamma)}.
\end{align*}
The term $\Vert p_h\Vert_{L^2(\Gamma)} \Vert \n\cdot \nabla p_h\Vert_{L^2(\Gamma)}$ can be controlled by $\Vert \nabla p_h\Vert^2_{L^2(\Omega)} $ and $h^{-1} \Vert p_h  \Vert^2_{L^2(\Gamma)}$ exactly in the same way as $\Vert \n \times \u_h \Vert_{L^2(\Gamma)} \Vert \nabla \times \u_h \Vert_{L^2(\Gamma)} $ is controlled by $\Vert \nabla \times \u_h \Vert^2_{L^2(\Omega)}$ and $h^{-1}\Vert \n \times \u_h \Vert^2_{L^2(\Gamma)} $ in Theorem~\ref{th:gal-nit}, now using the fact that $N_p$ is sufficiently large. This yields:
\begin{align*}
& B_{\rm SN} ([\u_h,p_h] , [\u_h,-p_h]) 
 \gtrsim \nu \Vert \nabla \times \u_h \Vert^2_{L^2(\Omega)} + \frac{L_0^2}{\nu} \Vert \nabla p_h\Vert^2_{L^2(\Omega)} + \frac{\nu h^2}{L_0^2 } \Vert \nabla\cdot\u_h\Vert^2_{L^2(\Omega)} \cl
& \qquad 
+ N'_u \frac{\nu}{h} \Vert \n \times \u_h \Vert^2_{L^2(\Gamma)} 
+ N'_p \frac{L_0^2}{\nu h} \Vert p_h  \Vert^2_{L^2(\Gamma)},
\end{align*}
for certain $N'_u\leq N_u$ and $N'_p \leq N_p$. 

Let now $[\vt_{h,0},q_{h,0}]\in V_{h,0}\times Q_{h,0}$ be the pair whose existence is established in Theorem~\ref{th-0-stab} that satisfies \Eq{gen-inf-sup-stab}, which we take such that  $\Vert [\vt_{h,0},q_{h,0}]\Vert_{V\times Q,{\rm S}} = \Vert [\u_{h,0},q_{h,0}]\Vert_{V\times Q,{\rm S}} \leq \Vert [\u_{h},p_{h}]\Vert_{V\times Q,{\rm S}}$. Using the fact that $\n\times \vt_{h,0} = {\bf 0}$ and $q_{h,0} = 0$ on $\Gamma$, now we get
\begin{align}
B_{\rm SN}([\u_h,p_h],& [\vt_{h,0} , q_{h,0}])
  = B_{\rm S}([\u_h,p_h],[\vt_{h,0} , q_{h,0}])
 - \nu \langle{\bm n}  \times \u_h , \nabla\times \vt_{h,0}\rangle_\Gamma \cl
&  - \langle {\bm n}  \cdot \vt_{h,0} , p_h\rangle_\Gamma
 - \frac{L_0^2}{\nu}\langle {\bm n}  \cdot \nabla q_{h,0} , p_h\rangle_\Gamma
 \cl
& \gtrsim \Vert [\u_h,p_h] \Vert^2_{V\times Q,{\rm S}}   
- \gamma \frac{\nu}{h}\Vert \n\times \u_h\Vert^2_{L^2(\Gamma)}
- \gamma \frac{L_0^2}{\nu h}\Vert p_h\Vert^2_{L^2(\Gamma)} \cl
& - \nu \langle{\bm n}  \times \u_h , \nabla\times \vt_{h,0}\rangle_\Gamma - \langle {\bm n}  \cdot \vt_{h,0} , p_h\rangle_\Gamma
 - \frac{L_0^2}{\nu}\langle {\bm n}  \cdot \nabla q_{h,0} , p_h\rangle_\Gamma. \el
\end{align}
The terms $\langle{\bm n}  \times \u_h , \nabla\times \vt_{h,0}\rangle_\Gamma $ and $\langle {\bm n}  \cdot \vt_{h,0} , p_h\rangle_\Gamma$ can be bounded as in Theorem~\ref{th:gal-nit}, just replacing the norm $\Vert [\u_{h},p_{h}]\Vert^2_{V\times Q}$ by $\Vert [\u_{h},p_{h}]\Vert^2_{V\times Q,{\rm S}}$, and the last term is also immediately bounded as 
\begin{align}
 \frac{L_0^2}{\nu}\langle {\bm n}  \cdot \nabla q_{h,0} , p_h\rangle_\Gamma
& \leq \frac{\alpha}{2} \Vert [\u_{h},p_{h}]\Vert^2_{V\times Q,{\rm S}}  + \frac{1}{2\alpha} {C^2_{\rm trace}} \frac{L_0^2}{\nu h} \Vert p_h \Vert^2_{L^2(\Gamma)},\el
\end{align}
for any $\alpha > 0$. The proof now proceeds in that of Theorem~\ref{th:gal-nit}.
\end{proof}

Once the inf-sup condition has been established, we may proceed to obtain stability and convergence. Let us start noting that the stabilisation terms do not modify the right-hand-side linear form, which is the same as for Nitsche's method using the Galerkin approach, i.e., the form $L_{\rm N}$ given by \Eq{gal-nit-L}. For this, we now have:

\begin{lemma}\label{lem:l3}
The linear form $L_{\rm N}$ given in \Eq{gal-nit-L} is continuous in the norm $\vtres{\cdot}_{V\times Q,{\rm SN}}$, the continuity constant being bounded as
\begin{align*}
\vtres{L_{\rm N}}_{{\mathcal L}(V_h\times Q_h,{\rm SN}; \mathbb{R})} \lesssim  \Vert {\bm f}\Vert_{V'} +  \Bigl(\frac{\nu}{h}\Bigr)^{1/2} \Vert \n\times \bar{\u}\Vert_{L^2(\Gamma)}.
\end{align*}
\end{lemma}

\begin{proof} It follows immediately from  $\vtres{[ \vt_h, q_h] }_{V\times Q,{\rm N}} \leq \vtres{[ \vt_h, q_h] }_{V\times Q,{\rm SN}}$
\end{proof}

\begin{lemma}\label{lem:l4}
For any $[\u,p]\in V\times Q$, let the interpolation error function be
\begin{align}
E_{\rm S}(\u,p;h) = \inf_{[ \tilde{\u}_h, \tilde{p}_h] \in V_{h}\times Q_{h}} D_{\rm S}([ \u - \tilde{\u}_h, p - \tilde{p}_h]),\label{eq:err-est-stab-nit}
\end{align}
where
\begin{align*}
D_{\rm S}([ \vt ,  q])
& = \vtres{ [\vt , q ] }_{V\times Q,{\rm SN}}  
+ {(\nu h)^{1/2}} \Vert \n\times \nabla \times  \vt \Vert_{L^2(\Gamma)} \cl
&+  \frac{(\nu h)^{1/2}}{L_0} \Vert \n\cdot \vt \Vert_{L^2(\Gamma)}
+ L_0 \Bigl(\frac{ h}{\nu}\Bigr)^{1/2} \Vert \n\cdot \nabla q \Vert_{L^2(\Gamma)}.
\end{align*}
Then, for all $[\vt_h,q_h]\in V_h\times Q_h$ there holds
\begin{align*}
\inf_{[ \tilde{\u}_h, \tilde{p}_h] \in V_{h}\times Q_{h}}  B_{\rm SN}( [\u - \tilde{\u}_h, p - \tilde{p}_h] , [\vt_h,q_h]) \lesssim E_{\rm S}(\u,p;h) \vtres{[\vt_h,q_h]}_{V\times Q,{\rm SN}}.
\end{align*}
\end{lemma}

\begin{proof} Following the proof of Lemma~\ref{lem:l1}, the only terms that deserve to be analysed in the expression of $B_{\rm SN}( [\u - \tilde{\u}_h, p - \tilde{p}_h] , [\vt_h,q_h])$ are:
\begin{align*}
\frac{L_0^2}{\nu} \langle \n\cdot \nabla q_h , p - \tilde{p}_h \rangle_\Gamma 
\lesssim \frac{L_0}{\nu^{1/2}} C_{\rm trace}  \Vert \nabla q_h \Vert_{L^2(\Omega)} \frac{L_0}{(\nu h)^{1/2}} \Vert p - \tilde{p}_h \Vert_{L^2(\Gamma)},  \cl
\frac{L_0^2}{\nu} \langle \n\cdot \nabla (p - \tilde{p}_h) , q_h\rangle_\Gamma 
\lesssim \frac{L_0 h^{1/2}}{\nu^{1/2}} \Vert \n\cdot \nabla (p - \tilde{p}_h) \Vert_{L^2(\Gamma)} \frac{L_0}{(\nu h)^{1/2}} \Vert q_h \Vert_{L^2(\Gamma)}.\el
\end{align*}
These terms are both bounded by $D_{\rm S}([\u -\tilde{\u}_h, p - \tilde{p}_h]) \vtres{[\vt_h,q_h]}_{V\times Q,{\rm SN}}$. 
\end{proof}

It is now immediate to show that $E_{\rm S}(\u,p;h)$ is the error function of the formulation:

\begin{theorem}\label{th:conv-stab-nit} 
Under the assumptions of Theorem~\ref{th:stab-nit}, problem \Eq{stab-nit} is well posed, in the sense that it admits a unique solution $[\u_h,p_h]\in V_{h}\times Q_{h}$ that satisfies
\begin{align}
\vtres{ [\u_h,p_h] \Vert}_{V\times Q,{\rm SN}} \lesssim \Vert {\bm f}\Vert_{V'} + \Bigl(\frac{\nu}{h}\Bigr)^{1/2} \Vert \n\times \bar{\u}\Vert_{L^2(\Gamma)}.\label{eq:gal-nit-stab}
\end{align}
Furthermore, $[\u_h,p_h]$ converges optimally as $h \to 0$ to the solution $[\u,p]\in V \times Q$ of the continuous problem \Eq{varc1}-\Eq{varc2}, in the following sense:
\begin{align}
& \vtres{ [\u - \u_h, p - p_h] }_{V\times Q,{\rm SN}} \lesssim  E_{\rm S}(\u,p; h),\el
\end{align}
where $E_{\rm S}(\u,p; h)$ is given in \Eq{err-est-stab-nit}.
\end{theorem}

\begin{proof} 
The same as that of Theorem~\ref{th:conv-gal-nit}.
\end{proof}

Using standard interpolation estimates, it is observed that the error functions of both the stabilised formulation, $E_{\rm S}(\u,p;h)$,  and the Galerkin formulation using inf-sup stable elements, $E(\u,p;h)$, have the same optimal asymptotic behaviour in terms of $h$, given by \Eq{asbeh}.

\section{Numerical examples}
\label{sec:numex}

In this section we provide some numerical results to confirm the theoretical findings on the convergence of Nitsche's method using the stabilised FE approximation given in Section \ref{subsec:nitschestabilised}. \rc We have chosen to test the stabilised formulation for two reasons. First, because the effect of Nitsche's method is the same as for the Galerkin method with inf-sup stable elements and, second, because Nitsche's method is particularly important when using continuous nodal based interpolations due to the conformity issue described in the Introduction. \ot

We consider approximating the solution to Problem \Eq{m1}-\Eq{m4} by means of the formulation given in \Eq{stab-nit-B} on three different domains--all in two dimensions. The method is applied with equal order of linear interpolations for all the unknowns on various types of triangular elements to be described below. In the simulations, the scaling coefficients that appear in \Eq{stab-nit-B} are taken as $N_u = N_p = 10^2$, for all the cases considered. The other characteristic values are given individually for each test in what follows. 

\subsection{The square domain} 
\label{subsec:sq}
\begin{figure}[!h]\setlength{\unitlength}{1cm}
\centering
\includegraphics[width=3.9cm, height=4.0cm]{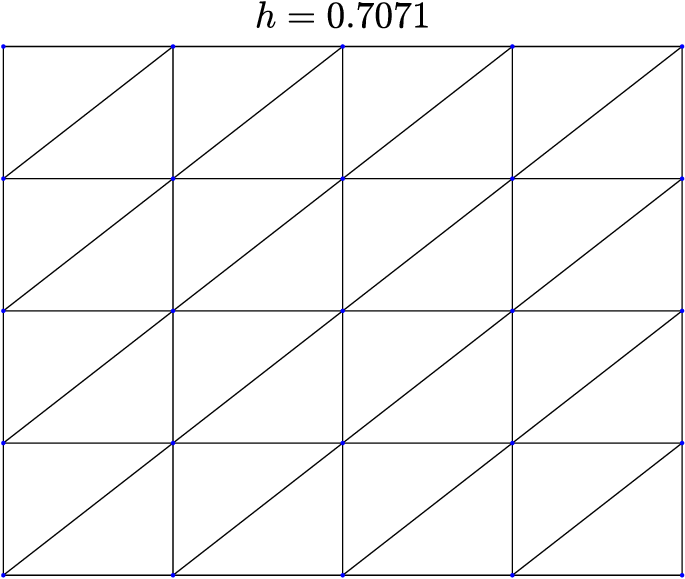}
\includegraphics[width=3.9cm, height=4.0cm]{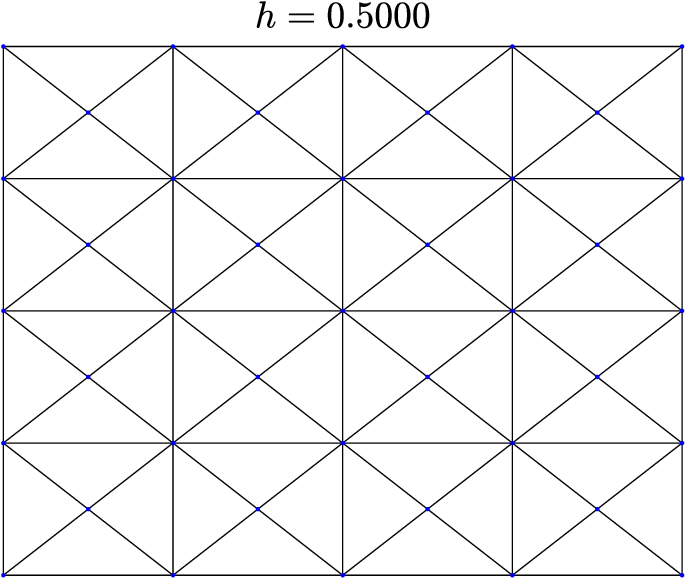}
\includegraphics[width=3.9cm, height=4.0cm]{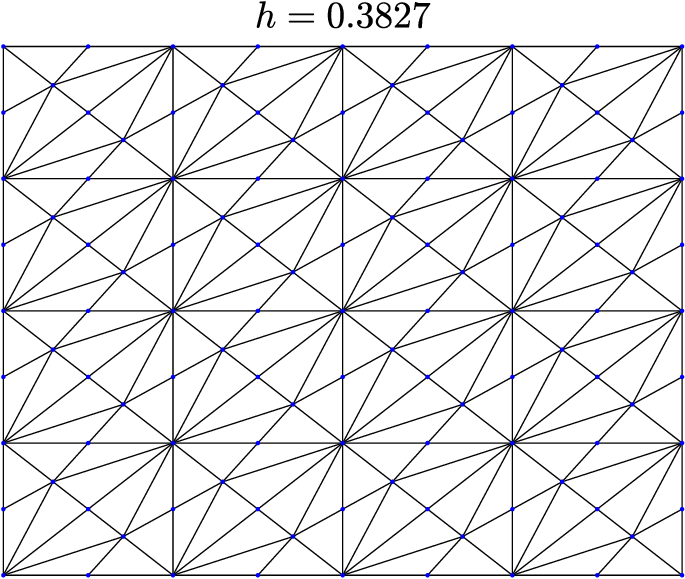}
\caption{The standard uniform right-angled (L), criss-cross (M), and Powell-Sabin (R) meshing of the square domain.}
\label{fig:sq_mesh}
\end{figure}

The first test problem is considered on the square domain $\Omega=(-1,1)^2$, with a smooth manufactured solution given by $\u(x,y) = (\varphi(x) \varphi'(y) , - \varphi'(x) \varphi(y))$, with $\varphi(t) = t^2 \sin(\pi t /2)$. This solution is used to determine $\bm f$, and then to check the convergence behaviour of the proposed scheme. We have performed the computations for this case on several mesh sequences, namely, standard uniform right-angled, criss-cross, and Powell-Sabin type meshes. Sample triangulations for these three mesh families are shown in Figure \ref{fig:sq_mesh}.
The characteristic length and the algorithmic stabilisation constant are taken respectively as $L_0=0.1$ and $c_u=0.1$, for the standard uniform right-angled mesh. The corresponding values are taken as $L_0=2$ and $c_u=1$, for the other two mesh sequences. We list the norms of the resulting numerical errors $\eu=\u - \u_h$ and $\ceu$ together with their rate of convergence towards zero as $h$ approaches zero in Table \ref{tab:sq_1}.

\begin{table}[!h]
\caption{Errors and rates of convergence (in brackets) for the square domain test on different triangulations.}  
\begin{center}
\begin{tabular}{c|c|ll}
\hline  Triangulation & $h$  & $\Vert \eu \Vert$  & $\Vert \ceu \Vert $
\tabularnewline  \hline  
Uniform right-angled &$ 0.3536 $ &   1.07e-01 & 9.64e-01 \\
&$ 0.1768  $&   2.04e-02 (2.39) & 4.31e-01 (1.16) \\ 
&$0.0884  $ &   4.75e-03 (2.10) & 2.15e-01 (1.00) \\ 
&$0.0442  $ &   1.18e-03 (2.00) & 1.08e-01 (1.00) \\ 
\hline
 Criss-cross &$ 0.2500 $ &  6.34e-02 & 3.91e-02 \\  
&$ 0.1250 $ &  1.60e-02 (1.98) & 1.00e-02 (1.96) \\   
&$ 0.0625 $ &  4.02e-03 (2.00) & 2.52e-03 (1.99) \\  
&$ 0.0312 $ &  1.01e-03 (2.00) & 6.31e-04 (2.00) \\   
\hline
Powell-Sabin &$ 0.1913 $ &   2.91e-02 & 2.63e-02 \\
&$ 0.0957 $ &   7.38e-03 (1.98) & 6.67e-03 (1.98) \\ 
&$ 0.0478 $ &    1.85e-03 (2.00) & 1.68e-03 (1.99) \\
&$ 0.0239 $ &    4.62e-04 (2.00) & 4.23e-04 (1.99) \\
\hline
\end{tabular}
\end{center}
\label{tab:sq_1}
\end{table}

It is evident from this table that the method is optimally convergent with double order of convergence in $\u_h$ for all the triangulations. The curl of the field also converges to  its expected value optimally for all the cases, while it exhibits a superconvergence in the case of special (criss-cross and Powell-Sabin) meshes. 
\begin{figure}[!h]\setlength{\unitlength}{1cm}
\centering
\includegraphics[width=4.8cm, height=4.2cm]{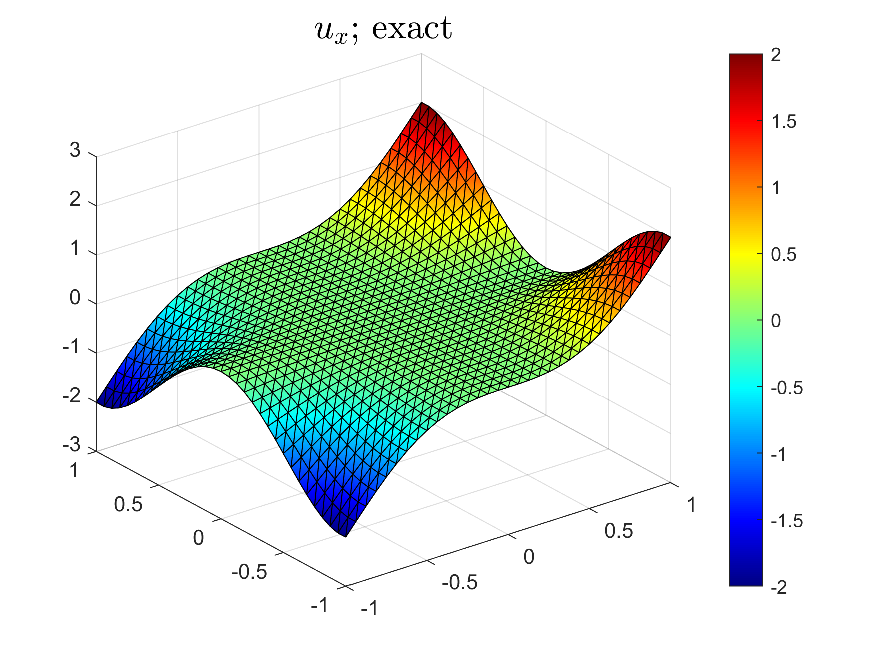}
\includegraphics[width=4.8cm, height=4.2cm]{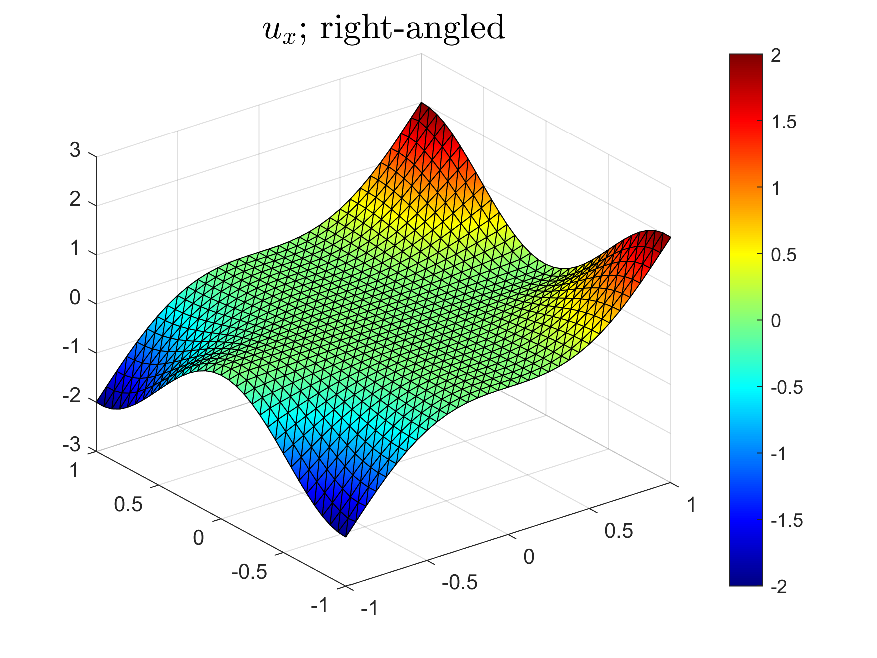}
\includegraphics[width=4.8cm, height=4.2cm]{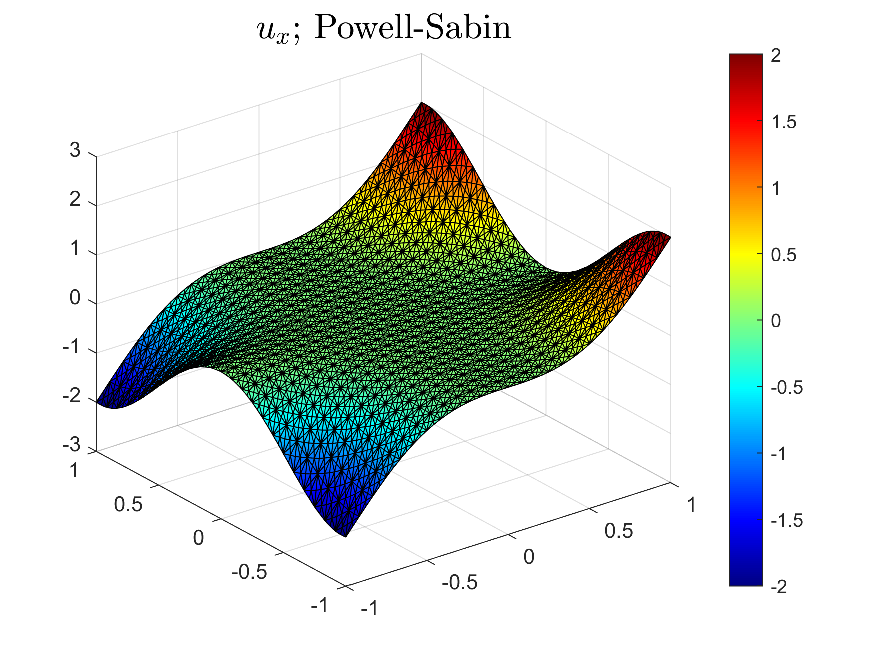}
\\
\includegraphics[width=4.8cm, height=4.2cm]{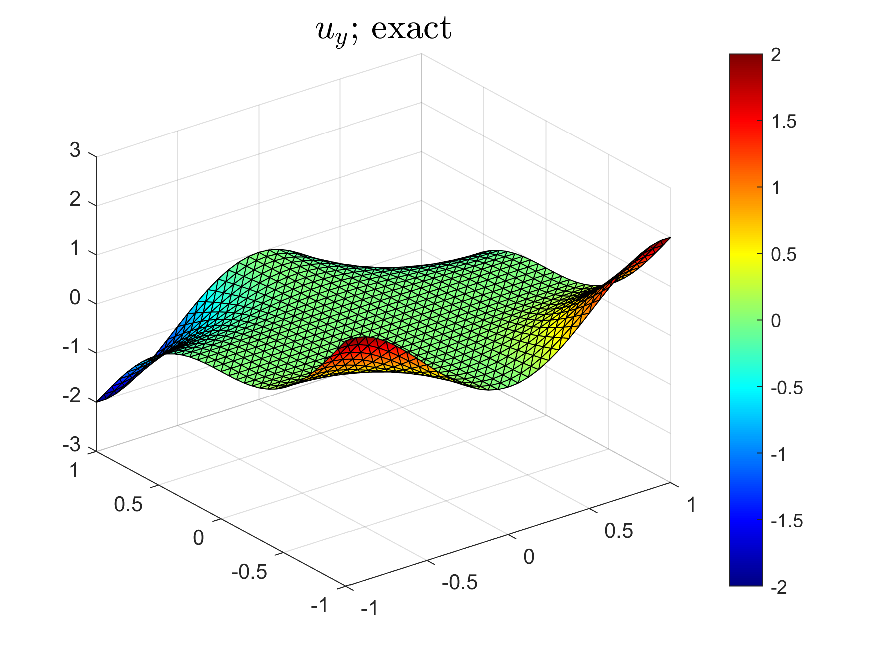}
\includegraphics[width=4.8cm, height=4.2cm]{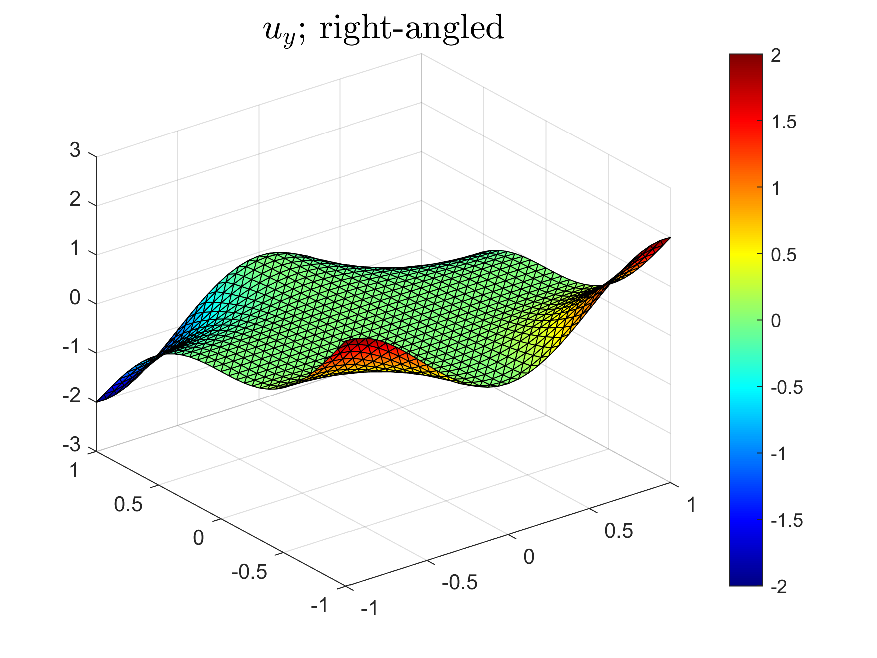}
\includegraphics[width=4.8cm, height=4.2cm]{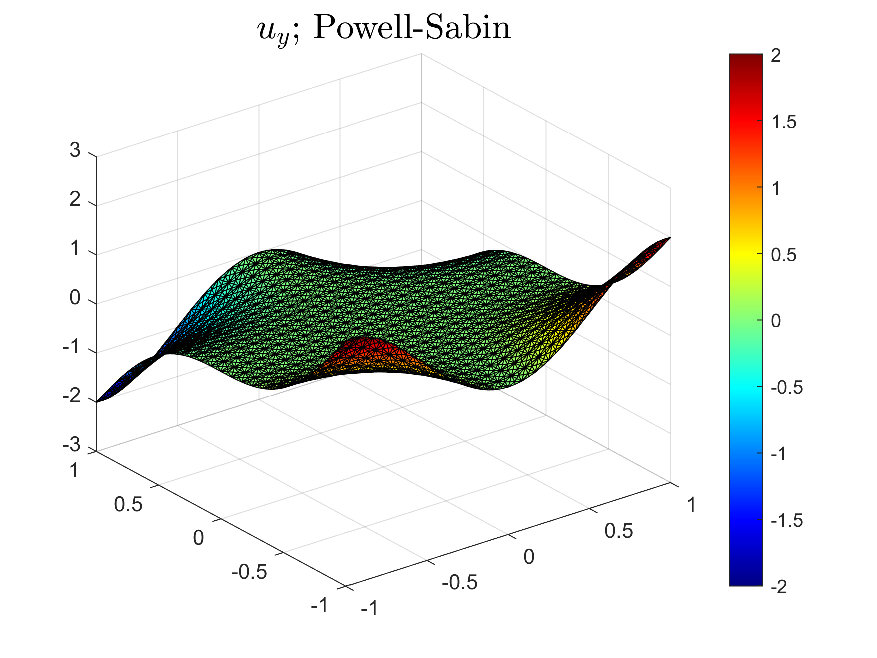}
\caption{Exact, right-angled mesh, and Powell-Sabin mesh solution components on the square domain.}
\label{fig:sq_components}
\end{figure}

To allow for a qualitative comparison of the computed solution components $u_x$ and $u_y$ with the exact ones, we present the surface plots of the exact (obtained on the right-angled mesh), right-angled mesh, and Powell-Sabin mesh solutions in Figure \ref{fig:sq_components}. The figure clearly shows a very good agreement between the computed and the analytical results. 

It is of significant importance in our study to check the comparison between the results obtained by weak prescription of the boundary conditions with those obtained by the strong imposition methodology. In order to do so,  we have considered the solution obtained on the Powell-Sabin mesh by strongly imposing the boundary conditions obtained making use of the exact solution. The results are listed in Table \ref{tab:sq_ps_strong}, and show clearly the close accuracy when compared with the corresponding ones obtained by Nitsche's method (given in Table~\ref{tab:sq_1}). 

\begin{table}[!h]
\caption{Strong imposition of the boundary conditions using  Powell-Sabin mesh on the square domain.}  
\begin{center}
\begin{tabular}{c|ll}
    \hline    $h$  & $\Vert \eu \Vert$  & $\Vert \ceu \Vert $
 \tabularnewline  \hline  
$ 0.1913 $ &   2.90e-02 & 2.65e-02 \\
$ 0.0957 $ &   7.38e-03 (1.98) & 6.73e-03 (1.98) \\ 
$ 0.0478 $ &   1.85e-03 (2.00) & 1.69e-03 (1.99) \\  
$ 0.0239 $ &   4.62e-04 (2.00) & 4.25e-04 (1.99) \\   
  \hline
\end{tabular}
\end{center}
\label{tab:sq_ps_strong}
\end{table}
 
\subsection{The L-shaped domain} 
\label{subsec:L}

In the second test, we consider a very widely used (e.g., in \cite{badia-codina-2009-2,houston-perugia-schotzau-2004-1}) configuration due to the presence of both smooth and nonsmooth solutions, the nonconvex domain defined by $\Omega=[-1,1]^2\setminus \{[0,1]\times[-1,0] \}$, with a re-entrant corner at the origin. The source function and the boundary conditions are taken so that the solution in polar coordinates is given as $\u=\nabla \psi$ where $\psi(r,\theta)=r^{2n/3} \sin(2n\theta/3)$, for different levels of smoothness depending on $n$. In our experiments, we consider the cases $n=1,2$, and $4$. For this example, we employ Nitsche's method using the stabilised formulation with $L_0=0.5$ and $c_u=1$. Due to the singularities involved, special types of meshes are necessary as we have already mentioned.  Thus, we use sequences of criss-cross and Powell-Sabin meshes to generate the results that are listed in Tables \ref{tab:L1cc} and \ref{tab:L1ps}, respectively.  

\begin{table}[!h]{\tiny
\caption{Errors and rates of convergence (in brackets) for the L-shaped domain test on criss-cross triangulations.}  
\begin{center}
 \begin{tabular}{c|ll|ll|ll}
\hline & \multicolumn{2}{c|}{$n=1$}  & \multicolumn{2}{c|}{$n=2$}   & \multicolumn{2}{c}{$n=4$}  \\  
    \hline    $h$  & $\Vert \eu \Vert$  & $\Vert \ceu \Vert $ & $\Vert \eu \Vert$  & $ \Vert \ceu \Vert$ & $\Vert \eu \Vert$    & $\Vert \ceu \Vert$
 \tabularnewline  \hline  
$0.1250$ & 2.61e-01 & 4.53e-01  & 2.12e-02 & 9.02e-02 & 3.09e-03 & 2.83e-02\\
$0.0625$ & 1.58e-01 (0.72) & 2.29e-01 (0.99)& 9.80e-03 (1.12) & 2.46e-02 (1.88)&8.33e-04 (1.89) & 3.68e-03 (2.94) \\
$0.0312$ & 9.38e-02 (0.76) & 1.02e-01 (1.17) &4.15e-03 (1.24) & 6.30e-03 (1.96)& 2.12e-04 (1.98) & 4.63e-04 (2.99)\\
$0.0156$ & 5.66e-02 (0.73) & 4.24e-02 (1.26) & 1.69e-03 (1.30) & 1.58e-03 (1.99)& 5.31e-05 (1.99) & 5.80e-05 (3.00)  \\
  \hline
\end{tabular}
\end{center}
\label{tab:L1cc}
}
\end{table}

\begin{table}[!h]{\tiny
\caption{Errors and rates of convergence (in brackets) for the L-shaped domain test on Powell-Sabin triangulations.}  
\begin{center}
\begin{tabular}{c|ll|ll|ll}
\hline & \multicolumn{2}{c|}{$n=1$}  & \multicolumn{2}{c|}{$n=2$}   & \multicolumn{2}{c}{$n=4$}  \\  
    \hline    $h$  & $\Vert \eu \Vert$  & $\Vert \ceu \Vert $ & $\Vert \eu \Vert$  & $ \Vert \ceu \Vert$ & $\Vert \eu \Vert$    & $\Vert \ceu \Vert$
 \tabularnewline  \hline  
$0.0957$ &   2.11e-01 & 3.48e-01    &   1.63e-02 & 5.03e-02    &  1.63e-03 & 1.34e-02   \\
$0.0478$ &  1.25e-01 (0.76) & 1.65e-01 (1.08)   & 6.94e-03 (1.23) & 1.09e-02 (2.20) & 4.27e-04 (1.93) & 1.71e-03 (2.97) \\
$0.0239$ &   7.40e-02 (0.76) & 7.09e-02 (1.22)  &2.81e-03 (1.30) & 2.23e-03 (2.30)   &1.08e-04 (1.99) & 2.12e-04 (3.01)\\
$0.0120$ &    4.49e-02 (0.72) & 2.91e-02 (1.29)   & 1.12e-03 (1.33) & 4.45e-04 (2.32)& 2.69e-05 (2.00) & 2.63e-05 (3.02) \\
  \hline
\end{tabular}
\end{center}
\label{tab:L1ps}
}

\end{table}It is clear from these tables that when $n=1$, the rate of convergence is determined by the regularity of the solution, as expected, since $\u \in H^{2n/3-\epsilon}(\Omega)^2$, for any $\epsilon >0$ for this problem \cite{badia-codina-2009-2,houston-perugia-schotzau-2004-1}. The same applies to the case $n=2$, in which it is still true that $\u \notin H^{2}(\Omega)^2$. On the other hand, when $n=4$, the solution $\u$ belongs to $H^{8/3-\epsilon}(\Omega)^2$, and with this smooth solution the error estimate applies optimally. All the results of these numerical investigations confirm the theoretical ones obtained in Section \ref{subsec:nitschestabilised}, and are in very good agreement with the associated ones reported in \cite{badia-codina-2009-2}.

Similar to what we have done in the previous example to compare the proposed weak prescription strategy with the strong imposition of boundary conditions, we intend to perform a final experiment for this case. However, the situation is more delicate for the present configuration due to the existence of a re-entrant corner and the utilisation of a nodal basis as we discussed earlier. Since we employ nodal interpolations, a number of alternatives can be considered to strongly impose the boundary condition \Eq{ma3}. A first option is to force both of the field components to be zero at the corner, and another option is to leave them free at this node. A third option can be achieved by defining a fictitious normal to the boundary, and adjusting the components so that the magnetic field follows the tangent to the boundary associated with this normal vector. For the critical case of $n=1$, we have implemented the described procedures and compare the resulting $L^2(\Omega)$ norms of the erro
 rs in Table \ref{tab:Lcomp}.

\begin{table}[!h]
\caption{A comparison of different ways to impose the boundary conditions on the L-shaped domain when $n=1$ and $h=0.0156$: Strong imposition (different strategies at the re-entrant corner) and Nitsche's method.}  
\begin{center}
\begin{tabular}{c|c|c}
 \hline      Strategy  & $\Vert \eu \Vert$  &    $\Vert \ceu \Vert $ 
\tabularnewline    \hline
$u_1=u_2=0$   & 5.66e-02  &  4.20e-02     \\
$u_1,u_2$ free & 2.82e-02  &  5.56e-03     \\
Bisector normal & 2.82e-02  &  5.56e-03     \\
\hdashline
Nitsche's method & 5.66e-02  &  4.24e-02    \\
\hline  
\end{tabular}
\end{center}
\label{tab:Lcomp}
\end{table}

The influence of different ways to prescribe the boundary condition on the numerical errors for this singular case can easily be observed from this table. As expected, the Nitsche method results are very close to the ones obtained by forcing both components to vanish at the origin. The other two strategies produce very similar results in terms of the computed errors.  
 
\subsection{The curved L-shape domain} 
\label{subsec:curvedL}

\begin{figure}[!h]\setlength{\unitlength}{1cm} 
\centering
\includegraphics[width=5.4cm, height=5.4cm]{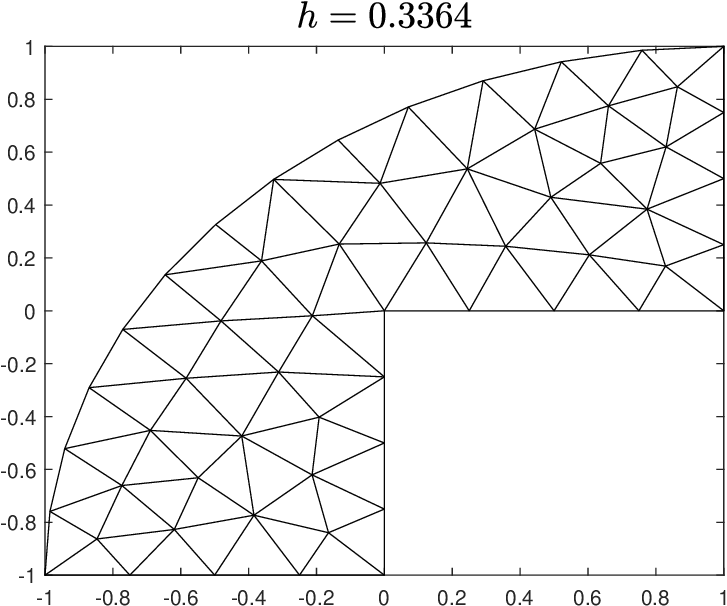}\quad
\includegraphics[width=5.4cm, height=5.4cm]{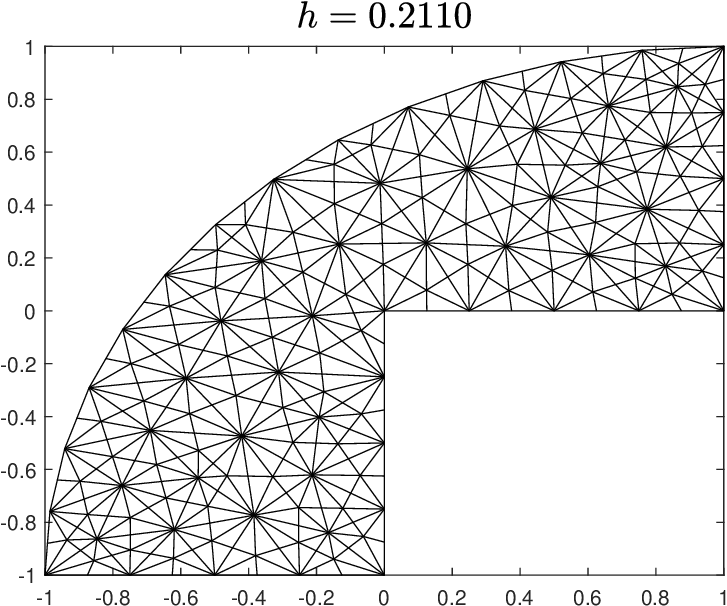}
\caption{The regular unstructured (L) and Powell-Sabin type (R) meshing of the curved L-shape domain.}
\label{fig:cL_mesh}
\end{figure}

As a last test, we consider the same solution as the previous L-shaped domain case on a curved L-shape domain now, and repeat the simulations whose results are presented in this subsection. The singularity occurring as a result of the re-entrant corner remains true as in the previous case. The significant difference in this one is the curved boundary that is obtained by joining the two diagonal corner points by a sector of a circle of radius 2, and centred at the point $(1,-1)$. The need for a weak prescription of boundary conditions is vital for this particular instance of a curved boundary. To discretise the computational domain, we have used a regular unstructured mesh and a sequence of Powell-Sabin type triangulations. Samples of both of these triangulations of the present domain are shown in Figure \ref{fig:cL_mesh}.

\begin{figure}[!h]\setlength{\unitlength}{1cm}
\centering
\includegraphics[width=4.8cm, height=4.2cm]{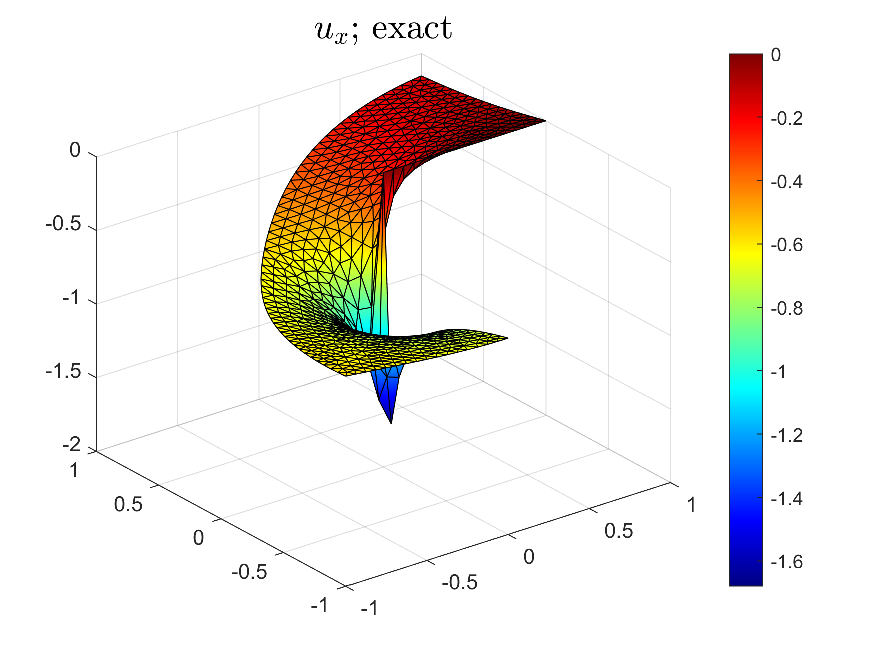}
\includegraphics[width=4.8cm, height=4.2cm]{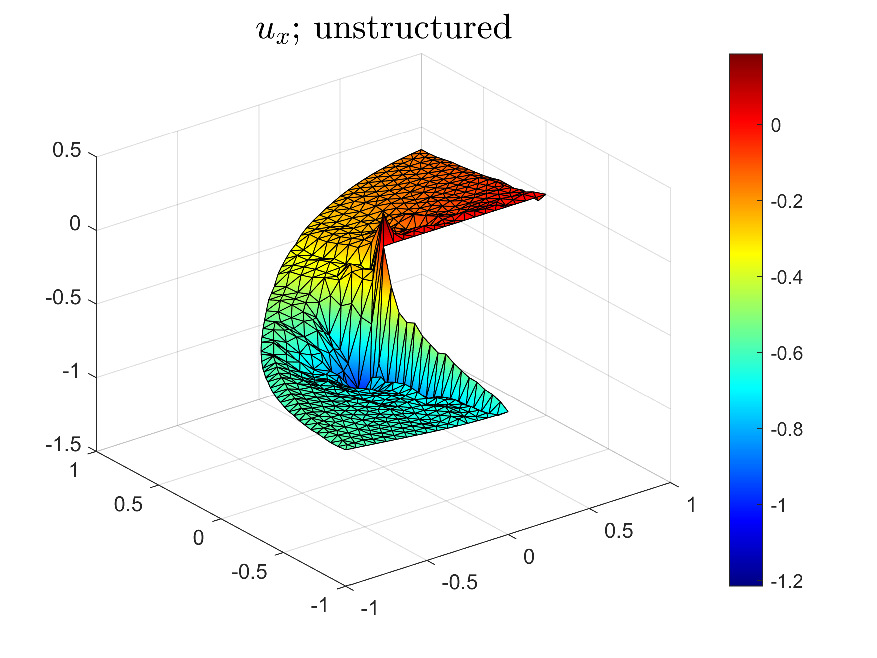}
\includegraphics[width=4.8cm, height=4.2cm]{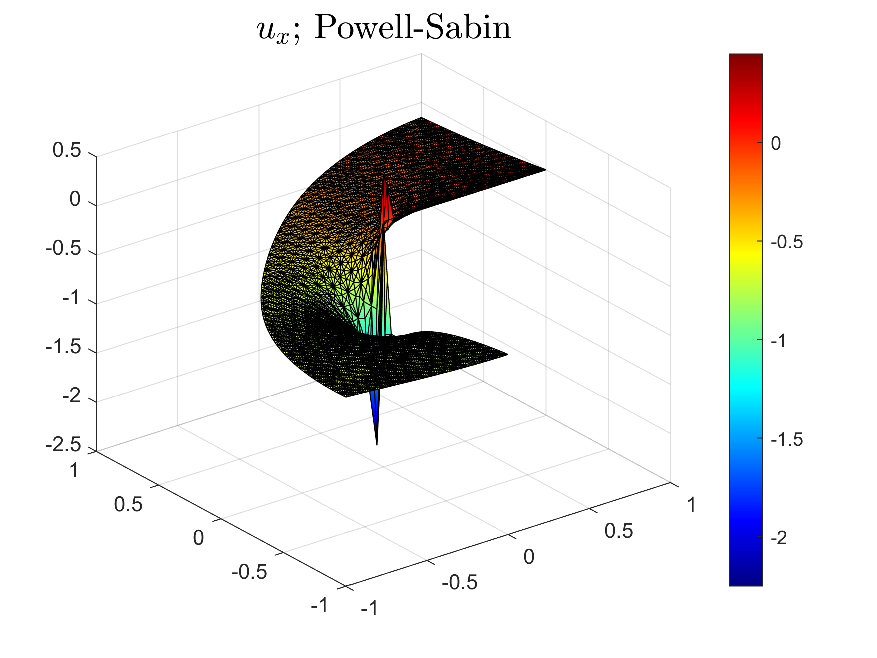}
\\
\includegraphics[width=4.8cm, height=4.2cm]{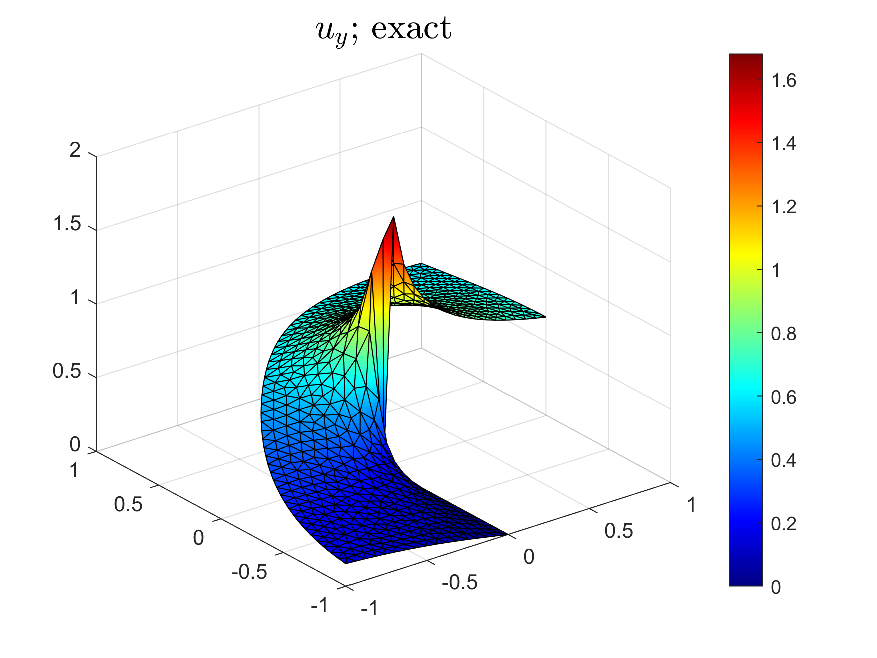}
\includegraphics[width=4.8cm, height=4.2cm]{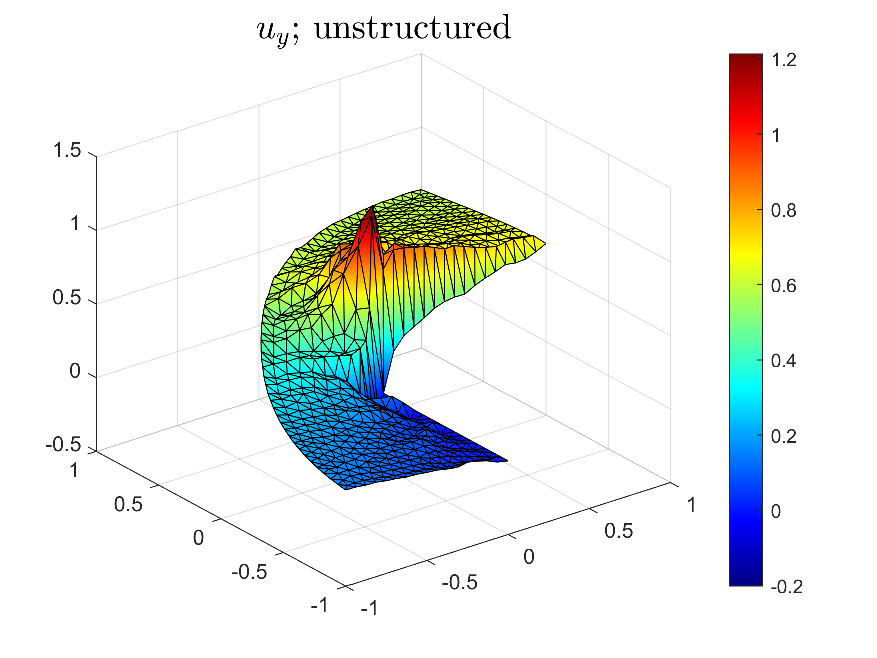}
\includegraphics[width=4.8cm, height=4.2cm]{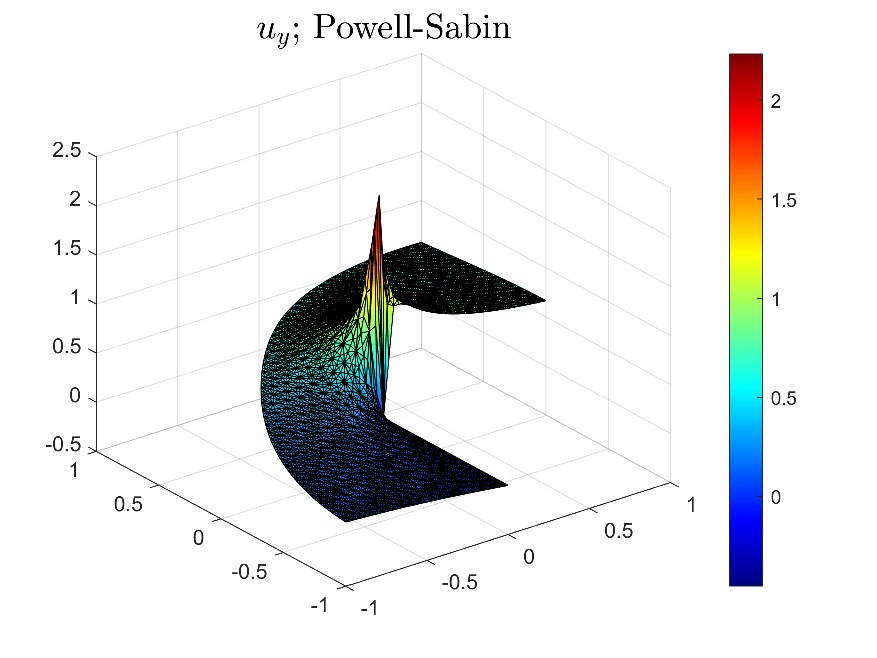}
\caption{Exact, unstructured mesh, and Powell-Sabin mesh solution components.}
\label{fig:cL_comp}
\end{figure}

We have carried out all the computations concerning this domain with the values $L_0=0.5$ and $c_u =0.1$. As already mentioned in Section \ref{sec:stab-exactBC}, the proposed scheme approximates the solutions with low Sobolev regularity optimally, provided that the used mesh has the ability to interpolate the corresponding scalar functions whose gradients are the solution components. Consequently, if this is not the case, the produced solutions may not capture accurately the correct solution behaviour. To explore this situation computationally, we have firstly used a standard regular unstructured mesh and then a Powell-Sabin type mesh to approximate the solution for the singular case when $n=1$. The results are presented in Figure \ref{fig:cL_comp} in terms of surface plots associated with the solution components of the magnetic field. The figure also depicts the corresponding exact solution generated on the regular unstructured mesh. 

The plots of the computed field components clearly demonstrate the failure of the employment of unstructured meshes in capturing the correct behaviour of the true solution. As we observed earlier, with the interpolation used on Powell-Sabin meshes, the expected solution is accurately recovered. 

In order to further analyse the convergence rates for this case, we have finally performed the computations on a sequence of these Powell-Sabin meshes, and list the errors with their rates in Table \ref{tab:cL_ps}.

\begin{table}[!h]{\tiny
\caption{Errors and rates of convergence (in brackets) for the curved L-shape domain test on Powell-Sabin triangulations.}  
\begin{center}
\begin{tabular}{c|ll|ll|ll}
\hline & \multicolumn{2}{c|}{$n=1$}  & \multicolumn{2}{c|}{$n=2$}   & \multicolumn{2}{c}{$n=4$}  \\  
    \hline    $h$  & $\Vert \eu \Vert$  & $\Vert \ceu \Vert $ & $\Vert \eu \Vert$  & $ \Vert \ceu \Vert$ & $\Vert \eu \Vert$    & $\Vert \ceu \Vert$
 \tabularnewline  \hline  
$0.4163$ & 3.15e-01 & 3.31e-01    &  4.49e-02 & 6.86e-02& 2.52e-02 & 2.00e-02 \\
$0.2110$ & 2.22e-01 (0.50) & 1.12e-01 (1.56)   & 2.66e-02 (0.76) & 2.02e-02 (1.77) & 7.79e-03 (1.69) & 3.08e-03 (2.70) \\
$0.1192$ & 1.42e-01 (0.65) & 4.31e-02 (1.38)   & 1.19e-02 (1.15) & 5.14e-03 (1.97) & 1.77e-03 (2.14) & 5.16e-04 (2.58) \\
$0.0586$ & 9.02e-02 (0.65) & 1.37e-02 (1.66)   & 5.33e-03 (1.16) & 1.09e-03 (2.24) & 4.90e-04 (1.85) & 9.96e-05 (2.37) \\
$0.0287$ & 5.59e-02 (0.69) & 4.20e-03 (1.70)   & 2.13e-03 (1.32) & 2.18e-04 (2.32) & 1.17e-04 (2.06) & 2.22e-05 (2.17) \\
$0.0146$ & 3.47e-02 (0.69) & 1.36e-03 (1.63)   & 8.07e-04 (1.40) & 4.51e-05 (2.27) & 2.82e-05 (2.06) & 5.44e-06 (2.03) \\
\hline
\end{tabular}
\end{center}
\label{tab:cL_ps}
}
\end{table}

These results once again show that the method attains the optimal convergence rates for all the different regularity levels considered, as anticipated from the theory. They also put forward that the convergence features are very similar to the ones obtained on the (straight) L-shaped domain.

\section{Conclusions}
\label{sec:conclu}

We have considered FE approximations of Maxwell's boundary value problem, and analysed the prescription of essential boundary conditions in a weak sense using Nitsche's method.  We have primarily focused on the analysis of two formulations with the inclusion of Nitsche terms; the Galerkin method when implemented with inf-sup stable elements, and an augmented-stabilised method that permits the use of nodal interpolations of arbitrary order. The analysis has been carried out by following a novel approach that relies on a splitting of the discrete spaces. We have provided the stability and convergence aspects for both of the formulations. 

In order to corroborate our theoretical findings in the case of the augmented-stabilised method, we have performed some numerical simulations. These results have confirmed the theoretical ones on optimal convergence of the method, and demonstrated the effectiveness of the proposed scheme in successfully approximating the expected solutions. In addition, the simulations have revealed the influence of the used meshes with different structures on correctly approximating the singular solutions. Finally, we have shown numerically that the results obtained by the weak prescription of the Dirichlet boundary conditions using Nitsche's approach comply well with the ones obtained by strong imposition.     
 
\bibliographystyle{plain}

\begin{thebibliography}{10}

\bibitem{arnold-1982-1}
D.N. Arnold.
\newblock An interior penalty finite element method with discontinuous
  elements.
\newblock {\em \Siamna}, 19:742--760, 1982.

\bibitem{assous1998}
F.~Assous, P.~Ciarlet, and E.~Sonnendr\"ucker.
\newblock Resolution of the {Maxwell} equations in a domain with reentrant
  corners.
\newblock {\em ESAIM: M2AN}, 32(3):359--389, 1998.

\bibitem{assous2011}
F.~Assous and M.~Michaeli.
\newblock Solving {Maxwell's} equations in singular domains with a {Nitsche}
  type method.
\newblock {\em Journal of Computational Physics}, 230(12):4922--4939, 2011.

\bibitem{aylwin2023sinum}
R.~Aylwin and C.~Jerez-Hanckes.
\newblock Finite-element domain approximation for {Maxwell} variational
  problems on curved domains.
\newblock {\em SIAM Journal on Numerical Analysis}, 61(3):1139--1171, 2023.

\bibitem{badia-codina-2009-2}
S.~Badia and R.~Codina.
\newblock A nodal-based finite element approximation of the {Maxwell} problem
  suitable for singular solutions.
\newblock {\em {SIAM} Journal on Numerical Analysis}, 50:398--417, 2012.

\bibitem{boffi2003}
D.~Boffi.
\newblock Finite elements for the time harmonic {M}axwell's equations.
\newblock In Peter Monk, Carsten Carstensen, Stefan Funken, Wolfgang Hackbusch,
  and Ronald H.~W. Hoppe, editors, {\em Computational Electromagnetics}, pages
  11--22, Berlin, Heidelberg, 2003. Springer Berlin Heidelberg.

\bibitem{bbf}
D.~Boffi, F.~Brezzi, and M.~Fortin.
\newblock {\em Mixed finite element methods and applications}, volume~44 of
  {\em Springer Series in Computational Mathematics}.
\newblock Springer, Heidelberg, 2013.

\bibitem{buffa-et-al-2002}
A.~Buffa, M.~Costabel, and D.~Sheen.
\newblock On traces for ${H}({\rm curl},{\Omega})$ in {Lipschitz} domains.
\newblock {\em Journal of Mathematical Analysis and Applications},
  276:845--867, 2002.

\bibitem{chaumont-frelet2023sinum}
T.~Chaumont-Frelet and P.~Vega.
\newblock Frequency-explicit a posteriori error estimates for finite element
  discretizations of {Maxwell's} equations.
\newblock {\em SIAM Journal on Numerical Analysis}, 60(4):1774--1798, 2022.

\bibitem{codina-baiges-2008-1}
R.~Codina and J.~Baiges.
\newblock Approximate imposition of boundary conditions in immersed boundary
  methods.
\newblock {\em \Ijnme}, 80:1379--1405, 2009.

\bibitem{codina-et-al-2009-1}
R.~Codina, J.~Principe, and J.~Baiges.
\newblock Subscales on the element boundaries in the variational two-scale
  finite element method.
\newblock {\em Computer Methods in Applied Mechanics and Engineering},
  198:838--852, 2009.

\bibitem{costabel-dauge-2002-1}
M.~Costabel and M.~Dauge.
\newblock Weighted regularization of {Maxwell} equations in polyhedral domains.
\newblock {\em \Nm}, 93:239--277, 2002.

\bibitem{O-douglas-dupont}
J.~Douglas and T.~Dupont.
\newblock Interior penalty procedures for elliptic and parabolic {Galerkin}
  methods.
\newblock In R.~Glowinski and J.~L. Lions, editors, {\em Computing Methods in
  Applied Sciences}, pages 207--216, Berlin, Heidelberg, 1976. Springer Berlin
  Heidelberg.

\bibitem{douglas-roberts}
J.~Douglas, Jr. and J.~E. Roberts.
\newblock Mixed finite element methods for second order elliptic problems.
\newblock {\em Mat. Apl. Comput.}, 1(1):91--103, 1982.

\bibitem{du2023}
Z.~Du and H.~Duan.
\newblock A stabilized finite element method on nonaffine grids for
  time-harmonic {Maxwell's} equations.
\newblock {\em BIT Numerical Mathematics}, 63:47, 2023.

\bibitem{hiptmair-2002}
R.~Hiptmair.
\newblock Finite elements in computational electromagnetism.
\newblock {\em Acta Numerica}, pages 237--339, 2002.

\bibitem{houston2005}
P.~Houston, I.~Perugia, A.~Schneebeli, and D.~Schötzau.
\newblock Interior penalty method for the indefinite time-harmonic {Maxwell}
  equations.
\newblock {\em Numerische Mathematik}, 100:485--518, 2005.

\bibitem{houston-perugia-schotzau-2004-1}
P.~Houston, I.~Perugia, and D.~Sch\"otzau.
\newblock Mixed discontinuous {Galerkin} approximation of the {Maxwell}
  operator.
\newblock {\em \Siamna}, 42:434--459, 2004.

\bibitem{juntunen-stenberg-09}
M.~Juntunen and R.~Stenberg.
\newblock Nitsche's method for general boundary conditions.
\newblock {\em Mathematics of Computation}, 78:1353--1374, 2009.

\bibitem{Kikuchi}
F.~Kikuchi.
\newblock Mixed and penalty formulations for finite element analysis of an
  eigenvalue problem in electromagnetism.
\newblock In {\em Proceedings of the first world congress on computational
  mechanics ({A}ustin, {T}ex., 1986)}, volume~64, pages 509--521, 1987.

\bibitem{monk-zhang-2020}
P.~Monk and Y.~Zhang.
\newblock {\em {\it Finite element methods for {Maxwell's} equations}, {\rm in
  Contemporary Mathematics 754: 75 Years of Mathematics of Computation}}.
\newblock American Mathematical Society, 2020.

\bibitem{perugias1}
I.~Perugia, D.~Sch\"otzau, and P.~Monk.
\newblock Stabilized interior penalty methods for the time-harmonic {Maxwell}
  equations.
\newblock {\em Computer Methods in Applied Mechanics and Engineering},
  191:4675--4697, 2002.

\bibitem{roppert2020}
K.~Roppert, S.~Schoder, F.~Toth, and M.~Kaltenbacher.
\newblock Non-conforming {Nitsche} interfaces for edge elements in
  curl--curl-type problems.
\newblock {\em IEEE Transactions on Magnetics}, 56(5):1--7, 2020.

\bibitem{stenberg-1995}
R.~Stenberg.
\newblock On some techniques for approximating boundary conditions in the
  finite element method.
\newblock {\em Journal of Computational and Applied Mathematics}, 63:237--339,
  1995.

\bibitem{wang2022}
N.~Wang and J.~Chen.
\newblock Convergence analysis of {Nitsche} extended finite element methods for
  {H}(curl)-elliptic interface problems.
\newblock {\em International Journal of Numerical Analysis and Modeling},
  19(4):487--510, 2022.

\bibitem{wheeler-1978}
M.~Wheeler.
\newblock An elliptic collocation-finite element method with interior
  penalties.
\newblock {\em SIAM Journal on Numerical Analysis}, 15:152--161, 1978.

\bibitem{winges2023}
J.~Winges and T.~Rylander.
\newblock Huygens' surface excitation for the finite element method applied to
  {Maxwell's} equations -- a construction based on {Nitsche's} method.
\newblock {\em Journal of Computational Physics}, 488:112237, 2023.

\end{thebibliography}

\end{document}